%% file: Deg2VF.tex
\documentclass[12pt,a4paper]{article}

\usepackage{amsmath,amsfonts,amssymb,enumerate}
\usepackage{amsthm,upref}
\usepackage{cancel}
\usepackage[a4paper]{geometry}
\usepackage{graphicx}
\usepackage{subcaption,tabularx,float}
\usepackage[final]{pdfpages}
\usepackage{tikz,tikzscale}
\usetikzlibrary{decorations.markings}
\usetikzlibrary{graphs,graphs.standard,arrows.meta}
\usetikzlibrary{arrows,shapes}
\usepackage[nodisplayskipstretch]{setspace}
\usepackage{pgfplots,pgfplotstable,pgf}
\usepgfplotslibrary{groupplots}
\pgfplotsset{compat=newest, ticks=none}
\usepackage{thmtools}
\usepackage{thm-restate}
\usepackage{etoolbox,environ}

\newcommand{\set}[1]{\left\lbrace #1 \right\rbrace}

\newcommand{\R}{\mathbb{R}}
\newcommand{\C}{\mathbb{C}}
\newcommand{\Q}{\mathbb{Q}}
\newcommand{\N}{\mathbb{N}}
\newcommand{\PP}{\mathbb{P}}
\newcommand{\nin}{\notin}
\newcommand{\pib}{\overline{\pi}}

\DeclareMathOperator{\Diff}{Diff_0}
\DeclareMathOperator{\Res}{Res}

\makeatletter
\newsavebox{\measure@tikzpicture}
\NewEnviron{scaletikzpicturetowidth}[1]{%
	\def\tikz@width{#1}%
	\def\tikzscale{1}\begin{lrbox}{\measure@tikzpicture}%
		\BODY
	\end{lrbox}%
	\pgfmathparse{#1/\wd\measure@tikzpicture}%
	\edef\tikzscale{\pgfmathresult}%
	\BODY
}
\makeatother

\newcounter{saveenumerate}
\makeatletter
\newcommand{\enumeratext}[1]{%
	\setcounter{saveenumerate}{\value{enum\romannumeral\the\@enumdepth}}
\end{enumerate}
#1
\begin{enumerate}
	\setcounter{enum\romannumeral\the\@enumdepth}{\value{saveenumerate}}%
}
\makeatother

\tikzset{cross/.style={cross out, draw=black, minimum size=2*(#1-\pgflinewidth), inner sep=0pt, outer sep=0pt},
	cross/.default={3pt}}

\theoremstyle{plain}
\newtheorem{thm}{Theorem}[section]
\newtheorem{proposition}[thm]{Proposition}
\newtheorem{corollary}[thm]{Corollary}
\newtheorem{lemma}[thm]{Lemma}

\theoremstyle{definition}
\newtheorem{definition}[thm]{Definition}
\newtheorem{example}[thm]{Example}

\tikzset{->-/.style={decoration={
			markings,
			mark=at position #1 with {\arrow{>}}},postaction={decorate}}}
		
\numberwithin{equation}{section}
\numberwithin{figure}{section}
\tikzset{every node/.append style={transform shape=false}}

\title{Regularisation for Planar Vector Fields}
\author{Nathan Duignan and Holger R.~Dullin}
\date{}

\begin{document}
	\maketitle
	\begin{abstract}
		This paper serves as a first foray on regularisation for planar vector fields. Motivated by singularities in celestial mechanics, the block regularisation of a generic class of degenerate singularities is studied. The paper is concerned with asymptotic properties of the transition map between a section before and after the singularity. Block regularisation is reviewed before topological and explicit conditions for the $ C^0 $-regularity of the map are given. Computation of the $ C^1 $-regularisation is reduced to summing residues of a rational function. It is shown that the transition map is in general only finitely differentiable and a method of computing the map is conveyed. In particular, a perturbation of a toy example derived from the 4-body problem is shown to be $ C^{4/3} $. The regularisation of all homogeneous quadratic vector fields is computed.
	\end{abstract}
	\let\clearpage\relax
	\include{Introduction}
	\include{Motivation}

	\include{Regularisation}

	\include{C0Regularisation}

	\include{RegularisingMap}
	\include{Normalising}

	\include{QuadraticVF}
	
	\appendix
	\include{Appendix}

	\bibliography{Deg2VF}
	\bibliographystyle{plain}
\end{document}

%% file: Introduction.tex
\section{Introduction}

Behaviour near singularities and the structure of their interconnections is the primary source of analysis for dynamical systems. Consequently, numerous techniques on qualitatively and quantitatively understanding local and global behaviour of singularities are known, see e.g.~\cite{andronov1974qualitative}. In this paper we describe and explore regularisation for singularities of planar dynamical systems. In essence, regularisation concerns the continuation of solutions to differential systems past singular points. Solutions that approach the singularity are called asymptotic orbits. If they can be extended through the singular point in some meaningful manner, then the singularity is deemed \textit{regularisable}. Historically, research of dynamical systems involved expanding solutions in power series of the independent variable $ t $. Extension of a solution past a singularity at $ t = t_s $ was naturally a question of whether, as a power series about $ t_s $, the solution had an analytic continuation. This is referred to as \textit{analytic regularisation}, but also appears in the literature as \textit{branch regularisation} and \textit{regularisation with respect to time}, see e.g. Sundman, Siegel \cite{siegel2012lectures} and more recently Wang and Punosevac \cite{punovsevac2012regularization} as well as Bakker et al \cite{bakker2011existence}.

Geometric analysis of dynamical systems seeks an understanding through the study of the flow of a vector field on a manifold. With this viewpoint, regularisation should be considered as a property of the flow, rather than of a particular solution. Conley and Easton \cite{conley1971isolated}, \cite{easton1971regularization} provide a precise definition of this notion, referred to as \textit{block regularisation} or \textit{regularisation by surgery}. Essentially, a singular set is deemed regular if there exists an extension of the singular orbits that is at least continuous with respect to initial conditions. This approach has been used most notably in celestial mechanics by Easton \cite{easton1972topology}, McGehee \cite{mcgehee1974triple}, Elbialy \cite{elbialy1990collision} and Simo and Martinez \cite{martinez1999simultaneous}. Whether the two notions are distinct is not clear from the definitions alone. However, McGehee \cite{mcgehee1974triple} has demonstrated that regularity in one sense need not imply regularity in the other.

Besides the general theory of dynamical systems and normal forms, see, e.g. \cite{Bruno1989},
the investigation of degenerate singularities of (planar) vector fields has at least two traditional roots. The older root is  
Hilbert's 16th problem on the number of limit cycles of polynomial vector fields, see, for instance,
\cite{Ilyashenko2008} and \cite{roussarie1995bifurcations}. Related to this root is work on the focus-centre problem, such as in \cite{manosa2002center}. 
The younger root is singularity theory where  holomorphic equivalence for complexified vector fields is studied.
For a review of this line of research see, e.g., \cite{VOBRG10,OBRGV14}, also see \cite{CG17}.
There is a classification theory in the complex, and first steps towards a related real 
classification have recently been made \cite{Jaurez-Rosas2017}. 
Our approach is different and less ambitious from both of these traditional roots,
in that we only consider a particular sub-class of singularities for which 
the notion of regularisability makes sense. As an interesting consequence we do find invariants of the real dynamics, in particular the 
degree of smoothness of the map ``across'' the singularity. It would be very interesting to bring together these  lines of investigation,
and identify where this finite degree of differentiability lies within the set of holomorphic invariants.

This paper is a first foray on the block regularisation of singularities of arbitrary planar vector fields. The primary aims are to review block regularisation of planar vector fields, classify which singular points of planar vector fields are regularisable, explore the differentiability of the resulting block map, and construct methods to compute the map. This is achieved for a generic class of vector fields. To the best of the Authors' knowledge, the study of asymptotic properties of the block regularisation of singularities has been restricted to specific singularities in celestial mechanics and not considered in a more general setting, such as planar vector fields. 

In Section \ref{sec:Motivation} a motivating toy example is derived from the simultaneous binary collision singularity in the 4-body problem. Then, block regularisation is reviewed in Section \ref{sec:BlockRegularisation} and Theorem \ref{thm:C0Regularisable} is proved, giving the first result on the $ C^0 $-regularisation of the planar vector fields. This theorem is substantiated in Section \ref{sec:HomogeneousVF} where the methods of blow-up and desingularisation are harnessed to prove Theorem \ref{thm:C0RegularisationHomX}. It is shown that $ C^0 $-regularity for a generic class of vector fields $ \mathcal{V} $ is determined by the number and sign of so called characteristic values. Section \ref{sec:CkRegularisation} begins the asymptotic analysis of the block map and constitutes the bulk of the paper. The block map is shown to be a composition of Dulac maps and smooth transition maps and hence, is in general a qausi-regular map. Theorem \ref{thm:NearId} and Corollary \ref{cor:Cauchy} provide results about the $ C^1 $-regularity and how it can be verified through the sum of the residues of a rational functional. This section also sees the construction of a method for direct computation of the block map. Finally, in Section \ref{sec:examples} a canonical form for homogeneous vector fields is proposed. It is used to classify the block regularisation of homogeneous quadratic vector fields. Ultimately, all of the tools developed are pulled together to prove that a perturbation of the motivating toy example is precisely $ C^{4/3} $-regularisable.

%% file: Motivation.tex
\section{Motivating Toy Example}\label{sec:Motivation}
	The motivation for this work comes from a curious conjecture about the nature of the simultaneous binary collision singularity in the 4-body problem \cite{martinez1999simultaneous}. Although it is resolved for the collinear problem \cite{martinez2000degree}, it remains to be shown for the planar problem. The conjecture states that any attempt at block regularisation of this singularity results in a block map that is precisely $ C^{8/3} $. This paper explores this odd loss of differentiability for planar systems. When does it occur? If so, how can it be computed?
	
\begin{figure}[ht]
		\centering
		\begin{scaletikzpicturetowidth}{\textwidth}
			\begin{tikzpicture}[scale = \tikzscale]
			\node[circle,fill,inner sep=2pt] (1) at (-5,0) {};
			\node[circle,fill,inner sep=1pt] (2) at (-2,0) {};
			\node[circle,fill,inner sep=2pt] (3) at (2,0) {};
			\node[circle,fill,inner sep=3pt] (4) at (4,0) {};
			\node[cross] (5) at (-4,0) {};
			\node[cross] (6) at (3,0) {}; 
			\node[inner sep=2pt,fill=white] (7) at (-4,0.5) {};
			\node[inner sep=2pt,fill=white] (8) at (3,0.5) {};
			
			\draw[-Latex] (1) -- node[label=south:$Q_1$] {}(2);
			\draw[-Latex] (3) -- node[label=south:$Q_2$] {}(4);
			\draw[-Latex] (7) -- node[auto] {$x$}(8);
			\end{tikzpicture}
		\end{scaletikzpicturetowidth}
		\caption{The configuration variables near simultaneous binary collision}\label{fig:3BDiff}
	\end{figure}
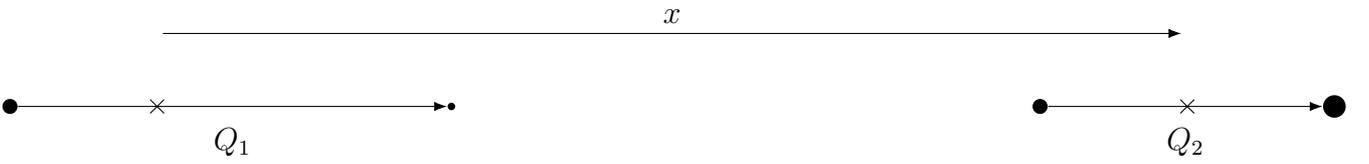

	Before proceeding to the precise definitions of regularisation and the analysis, it is useful to derive a toy example from the 4 body problem to help convey the theory. Suppose there are 4 collinear bodies consisting of two binaries undergoing collision in different regions of configuration space at precisely the same time $ t_f $. Further suppose that the bodies with mass $ m_1 $ and $ m_2 $ undergo one of the binary collisions and bodies with masses $ m_3,m_4 $ undergo the other. Let $ Q_1$ be the distance of $m_1$ and $m_2$ and  $Q_2 $  the distance between $m_3$ and $m_4$ and let $ x $ be the distance between the two centre of masses of the binaries, as seen in figure \ref{fig:3BDiff}. If $ P_1,P_2,y $ are the conjugate momenta of $ Q_1,Q_2,x $, the dynamics is determined by the differential equation,
	\begin{equation}
		\begin{aligned}
			\dot{Q}_1	&= M_1 P_1	\\
			\dot{Q}_2	&= M_2 P_2	\\
			\dot{x}		&= \mu y	\\
			\dot{P}_1	&= -k_1\frac{|Q_1|}{Q_1^3} + \frac{\partial K}{\partial Q_1}	\\
			\dot{P}_2	&= -k_2 \frac{|Q_2|}{Q_2^3} + \frac{\partial K}{\partial Q_2}	\\
			\dot{y}		&= \frac{\partial K}{\partial x}.
		\end{aligned}
	\end{equation}
	where the smooth function $ K = K(x,Q_1,Q_2) $ contains the potential terms coupling the two binaries and $ M_i,k_i,\mu > 0 $ are constant functions of the masses. 

	Note the singularity occurring when $ Q_1=0, Q_2=0 $ or $ Q_1=Q_2=0 $, the so called binary collisions and simultaneous binary collision, respectively. One can remove the singularity at 
	a binary collision through a Levi-Civita regularisation \cite{elbialy1993simultaneous}. This is done through the introduction of the generalised Levi-Civita variables $ (z_i,u_i) $ 
		\[ Q_i = \frac{1}{2}z_i^2,\qquad P_i = z_i^{-1} u_i, \]
	and a rescaling of time
		\[ d\tau = z_1^2 z_2^2 dt. \] 
	A second transformation of $ u_i $ to the intrinsic energy of the binary $ h_i $ is defined through
		\[ z_i^2 h_i = \frac{1}{2} M_i u_i^2 - 2 k_i. \]
	This series of transformations results in the vector field,
	\begin{equation}
		\begin{aligned}
			z_1^\prime	&= M_1 z_2^2 u_1(z_1,h_1)	\\
			z_2^\prime	&= M_2 z_1^2 u_2(z_2,h_2)	\\
			x^\prime	&= \mu z_1^2 z_2^2 y			\\
			y^\prime	&= z_1^2 z_2^2\frac{\partial K}{\partial x}\\
			h_1^\prime	&= M_1 u_1(z_1,h_1) z_2^2 \frac{\partial K}{\partial z_1}	\\
			h_2^\prime	&= M_2 u_2(z_2,h_2) z_1^2 \frac{\partial K}{\partial z_2},
		\end{aligned}	
	\end{equation}
	where the prime $\prime$ denotes the derivative with respect to the new time $\tau$.
	The binary collisions occurring when precisely one of the binaries is at collision ($ z_1= 0 $ or $z_2 = 0$) are now regular points of the vector field. However, the simultaneous binary collision $ z_1 = z_2 = 0 $ is transformed into an equilibrium point. Expansion about the simultaneous binary collision shows that at leading order in the small 
	quantities $z_1, z_2$ we have $ u_i(z_i,h_i) \sim 2 \sqrt{k_i/M_i} $ and $ \frac{\partial K}{\partial z_i} \sim D_i z_i $ for some $ D_i $  a constant depending on the masses. Looking at powers of  the small quantities $ z_1,z_2 $ one can see that the dynamics is primarily in the $ z_1,z_2 $ direction and hence, close to the collision, the dynamics is approximated by a foliation of 2-planes parallel to the $ z_1 z_2 $-plane. The dynamics on these 2-planes is the toy example.
	\begin{figure}[ht]
		\centering
		\includegraphics[width=0.4\textwidth]{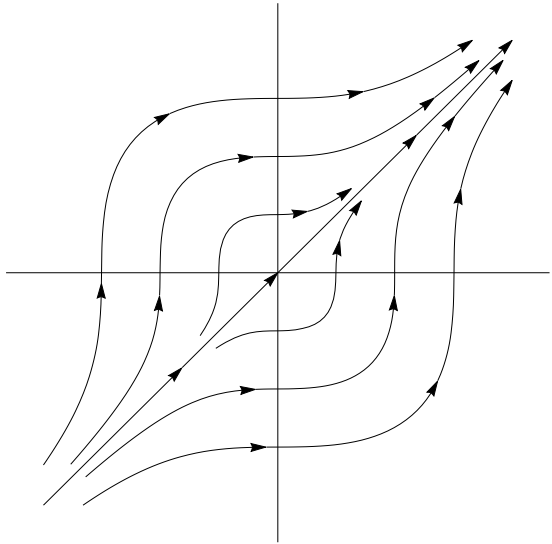}
		\caption{The phase portrait near the origin of the toy Example \ref{exam:Toy} for $ C_1=C_2=1 $.}\label{plot:toyexamplot}
	\end{figure}
\begin{example}[Toy Example]\label{exam:Toy}
		\begin{align*}
			z_1^\prime	&= C_1 z_2^2	\\
			z_2^\prime	&= C_2 z_1^2,			
		\end{align*}
		with $ C_1,C_2 >0 $.
	\end{example}
	
	A phase portrait is shown in Figure \ref{plot:toyexamplot}. Although the origin is an equilibrium point, the flow in a neighbourhood is monotonically increasing in the direction of the diagonal. One could hope to `pull out' the singularity and allow the asymptotic orbit to pass straight through, `regularising' the vector field.
This is made precise in the next section using the ideas of Conley and Easton.

%% file: Regularisation.tex
\section{Block Regularisation}\label{sec:BlockRegularisation}
	In this section block regularisation is presented. The key details and definitions from \cite{easton1971regularization} necessary for the proceeding analysis are given and we prove a theorem detailing necessary and sufficient conditions for $ C^0 $-block regularisation.
	
	Let $ V $ be a smooth vector field defined on a smooth manifold $ M $ and $ \varphi: M\times \R \to M $ the flow of $ V $. Note that it is assumed the flow is complete but the ideas presented extend to more general flows. A set $ I \subset M $ is called an invariant set if it is invariant under the flow of $ \varphi $, that is, $ \varphi(I,\R) = I $. Related to these invariant sets is the following definition.
	\begin{definition}
		An isolating neighbourhood for an invariant set $ I $ is an open set $ U \subset M $ containing $ I $ such that $ p \in \partial U $ implies $ \varphi(p,\R) \not\subset \overline{U} $. If  $ I $ is the maximal invariant set in an isolating neighbourhood $ U $ then we say $ I $ is isolated and $ U $ is an isolating neighbourhood of $ I $.
	\end{definition}
	Essentially, an isolating neighbourhood $ U $ of an invariant set $ I $ ensures the absence of another invariant set between $ I $ and $ \overline{U} $. Consequently, any open subset of $ U $ containing $ I $ must also be an isolating neighbourhood of $ I $.
	
	Let $ N $ be a smooth submanifold with boundary of $ M $ having $ \dim N = \dim M $ and let $ n := \partial N $, the boundary of $ N $. Then $ n $ can be partitioned into the following sets, see Figure \ref{fig:hyperbolicExample} for an example.
	\begin{definition}
		Define the following subsets of $ n = \partial N $.
		\begin{align*}
			n^+ 	&:= \set{p \in n\ |\ \exists t > 0 \text{ such that } \varphi(p,(-t,0))\cap N = \emptyset}\\ 
			n^- 	&:= \set{p \in n\ |\ \exists t > 0 \text{ such that } \varphi(p,(0,t))\cap N = \emptyset}\\ 
			\tau 	&:= \set{p \in n\ |\ V \text{ is tangent to $ n $ at $ p $}}\\
			a^+		&:= \set{p\in n^+\ |\ \varphi(p,t) \in N,\ \forall t \geq 0}\\
			a^-		&:= \set{p\in n^-\ |\ \varphi(p,t) \in N,\ \forall t \leq 0}			
		\end{align*}
	\end{definition} 
	$ n^+ $ contains all the points on $ n $ that flow into $ N $, $ n^- $ the points flowing out of $ N $, and $ \tau $ the points flowing tangent to $ N $. The sets $ a^+ $ and $ a^- $ contain the \textit{asymptotic points of $ N $}, which never flow out of $ N $ in forward and backward time respectively.
	
	\begin{definition}
		$ N $ is an \emph{isolating block} (for $ I $) if $ n^+ \cap n^- = \tau $, and if $ \tau $ is a smooth submanifold of $ n $ with codimension one and consequently $ n^+ $ and $ n^- $ are submanifolds of $ n $ with common boundary $ \tau $.
	\end{definition}

	By definition, an isolating block $ N $ is also an isolating neighbourhood of $ I $. The converse, that the existence of an isolating neighbourhood $ U $ implies the existence of an isolating block, is true \cite[p.~39]{conley1971isolated}.
	
	\begin{thm}[Conley and Easton]\label{thm:SetIffBlock}
		$ I $ is an isolated invariant set if and only if it is the maximal invariant set in some isolating block.
	\end{thm} 
	To prove the theorem, Conley and Easton show that for any isolated invariant set, transversal sections can be defined near the asymptotic sets $ a^+ $ and $ a^- $. Then, it is shown that these transversal sections to the flow can be extended using a collar to form an isolating block. Summarised in the following lemma is the specific requirements on the transversal section to ensure the existence of the collar.
	\begin{lemma}\label{lem:Sections}
		Let $ S_+ $ and $ S_- $ be disjoint local surfaces of section and let $ \tilde{n}^+ \subset S_+ $ be a compact submanifold with the same dimension as $ S_+ $ and with boundary $ \tilde{\tau} $. Assume also that $ \pi(\tilde{\tau}) $ (the projection of $ \tilde{\tau} $ onto $ S_- $) bounds a compact submanifold $ \tilde{n}^- \subset S_- $ and that $ T  $ is the orbit segment joining $ \tilde{\tau} $ and $ \pi(\tilde{\tau}) $. Finally, assume that $ \tilde{n} = \tilde{n}^+\cup T \cup \tilde{n}^- $ is a piecewise smooth manifold which bounds a submanifold $ \tilde{N} $ and that $ \tilde{n}^+ $ and $ \tilde{n}^- $ is the set of incoming and outgoing points of $ \tilde{N} $ respectively.
		
		Then there is an isolating block $ N $ containing $ \tilde{N} $ with the property that $ \tilde{n}^+ \subset n^+ $ and $ \tilde{n}^- \subset n^- $. Further, the asymptotic and invariant sets of $ \tilde{N} $ are the same as for $ N $.
	\end{lemma}
	
	Consider Figure \ref{fig:hyperbolicExample} which details the a hyperbolic saddle point $ \left( \alpha^2 x, -\beta^2 y \right), \alpha,\beta\in\R $ and a first example of an isolating block. Specifically, the disk is $ N $, the circle $ n $, $ \tau $ the points on $ n $ where the orbits are tangential, $ n^+ $ consists of the arcs transversal to the $ y $-axis, and $ n^- $ the arcs transversal to the $ x $-axis. The sets $ a^+ $ and $ a^- $ contain the intersection of $ n $ with the stable and unstable manifolds respectively.
	
	\begin{figure*}[ht]
		\centering
		\caption{}
		\begin{subfigure}[t]{0.5\textwidth}
			\centering
			\input{HyperbolicExample.tikz}
			\caption{A hyperbolic critical point isolated by the disk $ N $.}\label{fig:hyperbolicExample}
		\end{subfigure}%
		~ 
		\begin{subfigure}[t]{0.5\textwidth}
			\centering
			\input{TrivialExample.tikz}
			\caption{A trivial example of a $ C^\infty $-regularisable point.}\label{fig:RegularExample}
		\end{subfigure}
	\end{figure*}
	
	The key motivation for the study of isolating blocks is to obtain qualitative information about the flow inside the block through an understanding of the flow on the boundary. To achieve this correspondence we require the following definitions.
	
	\begin{definition}
		Let $ N $ be an isolating block. Define,
		\begin{align*}
			A^+ &:= \set{p\in N | \varphi(p,t)\in N \text{ for all }  t \geq 0}\\
			A^- &:= \set{p\in N | \varphi(p,t)\in N \text{ for all } t \leq 0}\\
			A	&:= A^+ \cup A^-,
		\end{align*}
		and note that $ a^+ = A^+\cap n^+, a^- = A^-\cap n^- $. In the simple example in Figure \ref{fig:hyperbolicExample} $ A^\pm $ are the stable/unstable invariant manifolds that are contained in the block.
	\end{definition}

	A useful result about the topology of these sets has been proved by Conley and Easton \cite[p.~42]{conley1971isolated}.
	\begin{lemma}\label{lem:deformationRetraction}
		$ n^+\setminus a^+, n^-\setminus a^- $ are both strong deformation retracts of $ N\setminus A $.
	\end{lemma}
	The proof uses the flow $ \varphi $ and the differentiability of the vector field to define a suitable deformation retraction. The lemma motivates the following definition.
	
	\begin{definition}
		A map $ \pi : n^+\setminus a^+ \to n^-\setminus a^- $ can be defined by flowing the point $ p \in n^+ $ to its image on $ n^- $. In fact, $ \pi $ is a homeomorphism; a result due to Conley \cite{conley1971isolated}.
	\end{definition} 
	
	The flow in Figure \ref{fig:hyperbolicExample} induces a map $ \pi $ that is smooth. However, no continuous extension to a map from all of $ n^+ $ to $ n^- $ exists. To see this, consider two points on the upper branch of $ n^+ $, arbitrarily close and on either side of $ a^+ $. The point on the left of $ a^+ $ will be mapped under $ \pi $ to the right branch of $ n^- $ (arbitrarily close to $ a^- $) whilst the point on the right of $ a^+ $ will be mapped under $ \pi $ to the right branch of $ n^- $ (arbitrarily close to the second $ n^- $). Hence, it is impossible to extend $ \pi $ to a continuous $ \pib $ as points on either side of $ a^+ $ will always get mapped to different branches.
	
	If instead we consider the degenerate vector field $\nobreak{ X = \left(x^2 + y^2,0 \right) }$ shown in Figure \ref{fig:RegularExample}, it can immediately be seen that the map $ \pi $ has a unique, smooth extension to the entirety of $ n^+ $ onto $ n^- $. Essentially, even though the origin of $ X $ is a singular point, the geometry of the phase portrait is precisely as though the origin were a regular point. The following definition of Easton \cite{easton1971regularization} clarifies this idea of regularity.

	\begin{definition}\label{def:regularisable}
		Suppose that $ N $ is an isolating block for $ I $. If $ \pi $ admits a unique $ C^k $ extension $ \pib : n^+ \to n^- $ then $ I $ is said to be $ C^k $\emph{-regularisable} (or has $ C^k $ regularity) and $ \pib $ is denoted the \emph{block map}.
	\end{definition}

	\noindent\textbf{Remark:} The extension $ \pib $ must be unique to eliminate cases that are topologically conjugate to a focus. Clearly these types of singularities should not be considered regular, yet, because all orbits in a neighbourhood of a focus are asymptotic orbits, the homeomorphism $ \pi $ is trivial. So, there could be many extensions to $ \pib $. Hence, the uniqueness rules out these types of singularities.

	In this paper we will see that many planar polynomial vector fields have a block map that is only finitely differentiable. To be more specific about the regularity, we have the following definition.
	
	\begin{definition}
		\begin{enumerate}
			\item A function $ f $ is in $ C^{0,\alpha}, \alpha\in[0,1] $ if \[ |f(x)-f(y)| \leq |x-y|^\alpha. \]
			\item A function $ f $ is in $ C^{k,\alpha}, \alpha \in [0,1], k \in \N $ if $ f^{(k)}(x) \in C^{0,\alpha} $.
			\item Define $ C^\alpha := C^{\lfloor \alpha \rfloor,\alpha-\lfloor\alpha\rfloor}, \alpha \in \R\setminus\N $.
			\item For Log-Lipschitz functions (e.g. $ x Log x $) the notation $ C^{0,\alpha} $ with undetermined $ \alpha $ is used.
		\end{enumerate}
	\end{definition}
	
	Whilst Easton's definition \ref{def:regularisable} of regularisable sets $ I $ works for arbitrary isolated singular sets, this paper is concerned only when $ I $ contains an isolated singular point (which for ease of notation will also be denoted $ I $) and only for planar vector fields. This restriction immediately provides a useful lemma.
	
	\begin{lemma}\label{lem:NonEmptyA}
		Let $ X $ be a planar vector field and suppose $ I $ is an isolated singular point of $ X $ with isolating block $ N $. Then $ a^+,a^- $ are non-empty.
	\end{lemma}
	\begin{proof}
		The lemma is shown by contradiction. Suppose $ a^+, a^- $ are empty. Then any $ p \in A^+ $ must remain inside $ N $ for all $t \in \R $ and hence $ p\in A^- $ as well. A symmetrical argument for any $ p \in A^- $ then guarantees that $ A^+ = A^- $. Moreover, since $ A^+, A^-$ are invariant sets and $ I $ is the maximal invariant set in $ N $, then $ I = A^+ \cup A^- =: A $. Using Lemma \ref{lem:deformationRetraction} it is concluded that $ n^+,n^- $ are strong deformation retracts of $ N \setminus I $. But this is impossible as $ N\setminus I  $ is homotopically equivalent to $ S^1 $ and $ n^+ $ is certainly not (it is a proper subset of a closed loop). This is a contradiction and, consequently, $ a^+ $ and $ a^- $ must be non-empty.
	\end{proof}
	
	This helpful lemma give rise to a proof of the first necessary and sufficient condition for an isolated singular point to be at least $ C^0 $-regularisable.
	\begin{thm}\label{thm:C0Regularisable}
		Let $ V $ be a $ C^1 $ planar vector field containing an isolated singular point $ I $. Then $ I $ is $ C^0 $-regularisable if and only if it has precisely one asymptotic orbit in each of $ A^+ $ and $ A^- $. 
	\end{thm}
	\begin{proof}
		Assume that $ I $ is $ C^0 $-regularisable. Then by definition there must exist an isolating block $ N $ with boundary $ n $ such that the homeomorphism $ \pi $ has a unique continuous extension $ \pib $. The uniqueness of $ \pib $ guarantees that $ a^+,a^- $ are compact submanifolds of $ n^+,n^- $ with codimension greater than or equal to one. As $ V $ is planar, then $ n $ is one dimensional and we conclude that $ a^+,a^- $ must be of dimension 0, that is, either a finite union of disjoint points or empty. However, Lemma \ref{lem:NonEmptyA} shows that as $ I $ is non-empty so too is $ a^+,a^- $.
		
		If $ a^+ $ is non-empty, we can choose a $ p^+ \in a^+ $ that is projected to some $ p^- = \pib(p^+) \in a^- $. Choose an open neighbourhood $ U^+ \subset n^+ $ of $ p^+ $ and define its projection to a neighbourhood of $ p^- $ by $ U^- = \pib(U^+) $. We demand that $ p^+ $ is the only point of $ a^+ $ contained in $ U^+ $ (and hence $ p^- $ the only point in $ a^- $ inside $ U^- $). This is justified by noting $ a^+ $ contains disjoint points. But then, $ U^+ $ and $ U^- $ are transversal sections with a continuous mapping $ \pib $ between them, thus Lemma \ref{lem:Sections} provides an isolating block $ \tilde{N} $ with asymptotic sets $ \tilde{a}^+,\tilde{a}^- $ containing only $ p^+,p^- $ respectively. As $ I $ must be the maximal invariant set in $ \tilde{N} $ we conclude that $ A^+, A^- $ contains only the orbits emanating from $ p^+ $ and $ p^- $ respectively.
		
		Conversely, assume that $ I $ is an isolated singular point with only two asymptotic orbits $ A^+ $ and $ A^- $. As $ I $ is isolated then Theorem \ref{thm:SetIffBlock} shows the existence of an isolating block $ N $ for $ I $. As there are only two asymptotic orbits then $ a^+ $ and $ a^- $ contain precisely one point each, say $ p^+ $ and $ p^- $. By defining $ \pib(p^+) = p^- $ we clearly obtain a unique, continuous extension of the map $ \pi $ of $ n^+ \setminus a^+ $ onto $ n^- \setminus a^- $. Hence, $ I $ is at least $ C^0 $ regularisable.
	\end{proof}

	\begin{example}[Toy]
		By studying the phase portrait of our toy example \ref{exam:Toy} we can see one orbit in each $ A^+ $ and $ A^- $. Theorem \ref{thm:C0Regularisable} then gives $ C^0 $-regularity of the block map as expected. However, whilst the figure hints that these are the only asymptotic orbits, it does not prove they are. In the proceeding section we will confirm this picture.
	\end{example}
	
	Theorem \ref{thm:C0Regularisable} immediately eliminates a large class of vector fields from consideration. Lemma \ref{lem:Sections} is also essential to the analysis in the remaining chapters. We can refocus the study of the block map from the entirety of the boundary of $ N $ to two sections of the flow intersecting the asymptotic orbits.

%% file: HyperbolicExample.tikz
\begin{tikzpicture}[scale = 1.5]
	\draw[->- = 0.25] (0,0) to (2,0);
	\draw[->- = 0.75] (0,2) to (0,0);
	
	\draw[->- = 0.25] (0,0) to (-2,0);
	\draw[->- = 0.75] (0,-2) to (0,0);

	\draw[domain=0.5:1.5,->- = 0.25] plot (\x,1/\x);
	\draw[domain=0.5:1.5,->- = 0.25] plot (\x,-1/\x);
	\draw[domain=-0.5:-1.5,->-=0.25] plot (\x,1/\x);
	\draw[domain=-0.5:-1.5,->-=0.25] plot (\x,-1/\x);
	
	\draw circle [blue, radius= sqrt(2)];
	
	\foreach \Point in {(1,1), (-1,1), (1,-1), (-1,-1)}{
		\node at \Point {\textbullet};
	}
	\foreach \Point in {(0,1.41421), (-1.41421,0), (0,-1.41421), (1.41421,0)}{
		\node at \Point {$ \circ $};
	}
	
	\node at (0.25,1.6) {$ a^+ $};
	\node at (1.6,0.25) {$ a^- $};
	\node at (-0.25,0.8) {$ A^+ $};
	\node at (0.8,-0.25) {$ A^- $};
	\node at (1-0.42,1+0.1) {$ n^+ $};
	\node at (1+0.1,1-0.42) {$ n^- $};
	\node at (1+0.2,1+0.2) {$ \tau $};

\end{tikzpicture}

%% file: TrivialExample.tikz
\begin{tikzpicture}[scale = 1.5]
	\draw (-2,0) to (2,0);
	\draw (0,2) to (0,-2);

	\foreach \k in {-1*sqrt(2),-0.75*sqrt(2),-0.5*sqrt(2),-0.25*sqrt(2),0,0.25*sqrt(2),0.5*sqrt(2),0.75*sqrt(2),1*sqrt(2)}{
		\draw[domain=-1.5:0,->- = 0.5] plot (\x,\k);
		\draw[domain=0:1.5,->- = 0.5] plot (\x,\k);
	}
	
	\draw circle [radius= sqrt(2)];

	\node at (0,1.4142) {\textbullet};
	\node at (0,-1.4142) {\textbullet};
\end{tikzpicture}

%% file: C0Regularisation.tex
\section{$ C^0 $ Regularisation}\label{sec:HomogeneousVF}
	
	A proof of a topological condition on the neighbourhood of a singular point to ensure $ C^0 $ regularity, Theorem \ref{thm:C0Regularisable}, was given in section \ref{sec:BlockRegularisation}. In this section we build the necessary tools for extension of this condition to an easily computable condition on a class of vector fields $ \mathcal{V} $. We approach the problem by first showing the result for a subclass of homogeneous vector fields $ \mathcal{V}_{hom} \subset \mathcal{V} $ . Then, a stability theorem is proved in section \ref{sec:stability} to extend the condition to the entirety of $ \mathcal{V} $.
	
	\subsection{Desingularisation}\label{ssec:desingularisation}
		Firstly, note that any hyperbolic singularity will not satisfy the topological conditions of Theorem \ref{thm:C0Regularisable}. Our search for regularisable singular points is restricted to degenerate singular points, i.e. points where the vector field and its Jacobian vanish. As such, we require tools to study these degenerate points.
		In particular we briefly describe the powerful desingularisation techniques presented in, for example \cite{Bruno1989}, \cite{Ilyashenko2008} and \cite{roussarie1995bifurcations}. The fundamental idea is to replace a degenerate singularity of the vector field $ X $ with a higher dimensional manifold. Assuming $ X $ is two-dimensional with singularity at the origin, the most intuitive way to achieve the desingularisation is to introduce the polar coordinate mapping $ \Phi:\mathbb{S}^1\times\R \to \R^2 $ given by $ \Phi(\theta,r) = (r \cos\theta,r\sin\theta) $. The pull-back $ X_\theta = \Phi^*X $ is a smooth vector field on a cylinder $ \mathbb{S}^1\times\R $ called the \textit{polar blow-up} of $ X $. Through this mapping, the point $ (0,0) $ is replaced by the circle $ \mathbb{S}^1 \times 0 $. In practice however, it is more useful to work in charts of $ \mathbb{S}^1\times\R $. This can be done through the \textit{directional blow-up},
		\begin{equation}
			\begin{aligned}
				x-\text{direction, } P_x^{-1}&:(\hat{x},\hat{y}) \mapsto (\hat{x},\hat{x}\hat{y}),\\
				y-\text{direction, } P_y^{-1}&:(\bar{x},\bar{y}) \mapsto (\bar{x} \bar{y}, \bar{y}).
			\end{aligned}
		\end{equation}
		
		The pull-back of these maps $ X_x = (P_x^{-1})^*X $ and $ X_y=(P_y^{-1})^*X $ describe $ X_\theta $ in charts. If $ X = (P(x,y),Q(x,y)) $ with $ P,Q $ homogeneous of degree $ n $, then a quick computation gives,
		\begin{equation}
			\begin{aligned}
				X_x	&= \hat{x}^{n-1} \left(\hat{x}P(1,\hat{y}), Q(1,\hat{y}) - \hat{y} P(1,\hat{y})\right) ,	\\
				X_y &= \bar{y}^{n-1} \left(P(\bar{x},1) - \bar{x}Q(\bar{x},1), \bar{y} Q(\bar{x},1)\right).
			\end{aligned} 
		\end{equation}
		The factors $ \hat{x}^{n-1} $ and $ \bar{y}^{n-1} $ show that the singularity at $ 0 $ has been mapped to a line of singularities for both $ X_x $ and $ X_y $, provided $ n > 1 $. We denote each of these lines the \textit{critical manifolds}. 
		
		One can \textit{desingularise} the vector fields by a time-scaling, hence considering 
		\[ \hat{X} := \frac{1}{\hat{x}^{n-1}} X_x,\qquad \bar{X} := \frac{1}{\bar{y}^{n-1}} X_y. \]
		The rescaling does not affect any topological analysis as we are only concerned with the phase portrait of the vector fields and not the orbits as functions of time. The hat and bar notation is used to remind the reader that, as will be shown in later sections, for $ C^0 $-regularisable singularities we can construct $ \hat{X} $ to contain a hyperbolic saddle on $ \hat{x} = 0 $ and $ \bar{X} $ to have no singularities on $ \bar{y} = 0 $, see Figure~\ref{fig:DirectionalBlowUp}.
		
		\begin{example}[Toy]\label{exam:desing}
			The toy example \ref{exam:Toy} is homogeneous of degree $ n=2 $. Computing the pull back of $ P_x^{-1} $ and $ P_y^{-1} $ and desingularising produces the vector fields
			\begin{equation}
				\begin{aligned}
					\hat{X} &= (C_1\hat{x}\hat{y}^2,C_2 - C_1\hat{y}^3), \\
					\bar{X}	&= (C_1 - C_2 \bar{x}^3,C_2 \bar{x}^2 \bar{y}) .
				\end{aligned}
			\end{equation}
			The singularity at $ (z_1,z_2) = 0 $ has been transformed to a line of singularities in both $ \hat{X} $ and $ \bar{X} $, namely, the critical manifolds $ \hat{x} = 0$, $\bar{y}=0 $, respectively. On these critical manifolds, the only singular point occurs at $ \hat{y}^* = (C_2/C_1)^{1/3} $ and its image on the $ y- $chart at $ \bar{x}^* = (C_1/C_2)^{1/3} $. Evaluating the eigenvalues at $ (0,\hat{y}^*) $ yields $ \lambda_1 = (C_1 C_2^2)^{1/3},\ \lambda_2= -3 (C_1 C_2^2)^{1/3} $ with corresponding eigenvectors $ v_1 = (1,0) $ and $ v_2=(0,1) $. Hence, the toy system has been desingularised into a critical manifold containing a hyperbolic point. Note that the one hyperbolic point corresponds to both the ingoing and outgoing asymptotic orbits. If one does the analysis using the polar blow-up instead, the ingoing and outgoing asymptotic orbits correspond to distinct hyperbolic points on the critical manifold (see Figure \ref{plot:polarblowup}).			
		\end{example}
		
		The crucial benefit of the blow-up is the resulting \textit{desingularisation} of the singular point. After an application of the blow-up the hope is to replace a degenerate singularity by a critical manifold containing elementary singularities \cite{Dumortier1977}. An isolated singularity is \textit{elementary} if it is either hyperbolic or semi-hyperbolic, that is, the singularity has at least one eigenvalue with non-zero real part.
		
 However, after just one application of the blow-up, one may still encounter non-elementary singularities on the critical manifold. But, one can then apply the procedure again to these non-elementary singularities. Thankfully, Dummortier \cite{Dumortier1977} has shown that, under reasonable non-degeneracy conditions, a singular point in the plane requires only finitely many blow-ups in the desingularisation.
		
%		This paper aims to be a first foray on regularisation for planar vector fields and not a complete generalisation. 
We seek to minimise the complexity in the computation of the regularising map for conciseness and clarity. Hence, we consider the class of vector fields in which precisely one blow-up for desingularisation to only hyperbolic singular points on the critical manifold are considered. In the next section, Proposition \ref{prop:oneblowup} shows this class is equivalent to $ \mathcal{V}^{hom} $.
		
	\subsection{Conditions for $ C^0 $-regularisation}\label{sec:conditionsC0}
	
		In this section, Theorem \ref{thm:C0Regularisable} is used to give easily computable neccessary and sufficient conditions for $ C^0 $-regularity of vector fields in $ \mathcal{V}^{hom} $ (see definition \ref{def:Vhom}). The primary challenge is determining how many asymptotic orbits a singularity has. The problem has many similarities to the focus-centre problem, studied by many authors, for example in \cite{manosa2002center}. The following definitions arise from the work on the focus-centre problem. 
		
		\begin{definition}
			Let $ \gamma(t) $ be an orbit of $ X $ that tends to the origin as $ t \to \infty $ (or $ -\infty $) such that
			\[ \lim_{t\to\infty (-\infty)} \frac{\gamma(t)}{|\gamma(t)|} \in \mathbb{S}^1. \]
			Then $ \gamma(t) $ is called a \textit{characteristic orbit}. 
		\end{definition}
		In Section \ref{sec:BlockRegularisation} any orbit $ \gamma(t) $ approaching a singularity was defined an \textit{asymptotic orbit}. Note that all characteristic orbits are asymptotic orbits but the converse does not hold. For example, any linear focus will not have a limit in $ \mathbb{S}^1 $. Thus, despite all orbits asymptotically approaching the singularity, they are not characteristic orbits.
		
		Let $ X = (P(x,y), Q(x,y))$ be a homogeneous vector field. By studying the limiting behaviour of the ratio $ y(t)/x(t) $ as $ t \to \infty $, potential characteristic orbits can be found; 
		\[ \frac{d}{dt}\left(\frac{y}{x}\right) = \frac{x Q(x,y) - y P(x,y)}{x^2}. \]
		Hence, any root $ \omega $ in projective space $ \PP $ of the polynomial $ x Q(x,y) - y P(x,y) $ represents the direction of an invariant line and given the appropriate name \textit{characteristic direction}. If $ \omega $ is a root of multiplicity $ k $, then we say it is a \textit{characteristic direction of multiplicity $ k $}. 
		
		In order to find characteristic directions, it is most natural to introduce polar coordinates $ (r \cos(\theta),r \sin(\theta)) $. This reduces the directional polynomial to a function $ F(\theta) $, the zeroes of which represent the angle of a characteristic direction. However, $ F(\theta) $ is a degree $ n+1 $ polynomial in trig functions, making the study of characteristic directions unnecessarily complex. In practice, it is easier to deal with directional coordinates of $ \PP $,
		\begin{equation}\label{eqn:dcoords}
			(x,y) \mapsto (x,y/x) =: (\hat{x},\hat{y})
		\end{equation}
		Notice that this transformation is precisely $ X_x $, the $ x $-directional blowup of $ X $. The directional polynomial under $ X_x $ is a real polynomial of degree $ n+1 $ of the form 
		\begin{equation}
			p(\hat{y}):= Q(1,\hat{y}) - \hat{y} P(1,\hat{y}).
		\end{equation}
		There may be a root of $ p(\hat{y}) $ at $ \infty $. One can use the other chart $ (x/y,y) $ to view this root, or simply rotate the vector field $ X $. 
		
		The following known result, proved in \cite{andronov1974qualitative}, gives the relation between characteristic orbits and characteristic directions.
		\begin{proposition}[\protect{\cite{andronov1974qualitative}}]\label{prop:asympOrbitImpliesDirection}
			Let $ \gamma(t) $ be a characteristic orbit of the origin and $ \displaystyle \omega = \lim_{t\to\infty} \frac{\gamma(t)}{|\gamma(t)|}. $ Then $ \omega $ is a characteristic direction for $ X $.
		\end{proposition}
		
		Interestingly, the converse does not hold; see \cite{manosa2002center} for a counterexample. A converse statement can be found with restrictions on the class of vector fields and multiplicity of the characteristic direction.
		\begin{definition}\label{def:Vhom}
			Call a singular point of a vector field $ X $ \textit{isolated} if there are no other singular points in a neighbourhood of the point. Define the \textit{class of multiplicity one vector fields} $ \mathcal{V}^{hom} $ containing all homogeneous vector fields $ X $ with an isolated singular point that has only characteristic directions of multiplicity one.
		\end{definition}
		\begin{proposition}\label{prop:DirectionImpliesOrbit}
			Assume that $ X \in \mathcal{V}^{hom} $ has an isolated singular point at the origin. Then, there exists a characteristic orbit with $ \displaystyle \omega = \lim_{t\to\infty} \frac{\gamma(t)}{|\gamma(t)|} $ if and only if there is a characteristic direction $ \omega $.
		\end{proposition}
		\begin{proof}
			The forward implication follows directly from Proposition \ref{prop:asympOrbitImpliesDirection}.
			
			Conversely, if $ \omega $ is a characteristic direction of multiplicity one, then by definition there must exist some root $ \hat{y}^* $ of $ p(\hat{y}) $ of multiplicity one. By rotating $ X $, we can assume without loss of generality that $ \hat{y}^*=0 $ and hence $ p(\hat{y}) = \hat{y}\bar{p}(\hat{y}) $ for some degree $ n $ polynomial $ \bar{p} $. Take the desingularised system,
			\begin{equation}\label{eqn:transys}
				\begin{aligned}
					\hat{x}^\prime	&=  \hat{x} P(1,\hat{y})\\
					\hat{y}^\prime	&=  \hat{y}\bar{p}(\hat{y}),
				\end{aligned}
			\end{equation}
			The multiplicity of $ \hat{y}^* $ implies $ \bar{p}(0) \neq 0 $. Further, it must be that $ P(1,0) \neq 0 $, else the line $ y = 0 $ would be a line of singular points, contradicting that $ (0,0) $ is an isolated singular point of $ X $. With these two assumptions, it can be concluded that $ (0,0) $ is a hyperbolic singular point of \eqref{eqn:transys}. Hence, there exists an orbit distinct from $ u = 0 $ that approaches the origin as $ t \to \infty $ or $ -\infty $. By inverting the transformation \eqref{eqn:dcoords}, it is verified that this orbit corresponds to a characteristic orbit of $ X $ with characteristic direction $ \omega $.
		\end{proof}
		
		The assumption that $ X\in\mathcal{V}^{hom} $ is not a necessary condition, however, vector fields in $ \mathcal{V}^{hom} $ are precisely those with the desired nice desingularisation properties discussed at the end of Section \ref{ssec:desingularisation}.
		
		\begin{proposition}\label{prop:oneblowup}
			The class of homogeneous vector fields $ \mathcal{V}^{hom} $ is equivalent to the class of homogeneous vector fields that require precisely one blow-up for desingularisation to either hyperbolic singular points on the critical manifold or no singular points.
		\end{proposition}
		\begin{proof}
			The critical manifold is given by $ \hat{x}=0 $. We need to show that any singular point on this manifold is hyperbolic. From the proof of Proposition \ref{prop:DirectionImpliesOrbit} we can assume that any singular point is at $ \hat{y} = 0 $ and the desingularised vector field takes the form,
			\[ \hat{X} = \left( \hat{x} P(1,\hat{y}), \hat{y}\bar{p}(\hat{y}) \right). \]
			From the same proof, if $ X \in \mathcal{V}^{hom} $ then $ \bar{p}(0) \neq 0 $ and $ P(1,0) \neq 0 $. Hence, the singular point at $ \hat{y} = 0 $ on the critical manifold $ \hat{x} $ must be hyperbolic.
		\end{proof}
	
		We remark that for any $ \phi\in\Diff $, the group of diffeomorphisms fixing $ 0 $ and mapping a neighbourhood of the origin to itself, $ \phi^* X \in \mathcal{V}^{hom} $. The number of characteristic directions $ \hat{y}^* $ of $ X $ is also preserved under transformations $ \phi $. Lastly, for each $ \hat{y}^* $, if $ \lambda_1 $ is the eigenvalue corresponding to the direction along the critical manifold $ \hat{x} = 0 $ and $ \lambda_2 $ the other eigenvalue, then the \textit{ratio of hyperbolicity} $ \lambda_1/\lambda_2 $ is preserved. This motivates the following definition:
		
		\begin{definition}\label{def:criticalvalues}
			Let $ {\hat{y}^*_1,\hat{y}^*_2,\dots,\hat{y}^*_d} $ be the critical directions of $ X \in \mathcal{V} $ and let $ \lambda_{i,1} $ be the eigenvalue at $ (0,\hat{y}^*_i) $ corresponding to the direction along the critical manifold $ \hat{x} = 0 $ and $ \lambda_{i,2} $ the other eigenvalue. Let $ r_i = \lambda_{i,1}/\lambda_{i,2} $ be the \textit{ratio of hyperbolicity}. Then, $ r(X) = \set{r_1,r_2,\dots,r_d} $ are called the \textit{critical values} of $ X $.
		\end{definition}
		Consider the set of equivalence classes $ [X] $ with $ X \sim X' \iff X = K\phi^*X' $ for some $ \phi\in\Diff $ and $ K $ a smooth function. A natural question arises from Definition \ref{def:criticalvalues}; can one assign a unique set of critical values to each equivalence class $ [X] $? This would produce a result in the same vein as a Jordan normal form. The question is addressed more directly in Section \ref{sec:examples}. It is addressed in the literature on holomorphic vector fields on $ \C^2 $, see for instance \cite{Jaurez-Rosas2017}. In this literature $ r(X) $ is often referred to as the Camacho-Sad index.
		
		We are ready to prove the main theorem of this section.
		\begin{thm}\label{thm:C0RegularisationHomX}
			Let $ X \in \mathcal{V}^{hom} $. Then $ X $ is $ C^0 $-regularisable if and only if $ r(X)=\set{r^*} $ and $ r^* <0 $.
		\end{thm}
		\begin{proof}
			The theorem follows from Theorem \ref{thm:C0Regularisable} in which an isolated singular point was shown to be regularisable if and only if it has precisely one ingoing and one outgoing asymptotic orbit. Firstly assume that $ X $ is $ C^0 $-regularisable. Proposition \ref{prop:asympOrbitImpliesDirection} guarantees the existence of a unique characteristic direction $ \hat{y}^* $ corresponding to both of the asymptotic orbits. By rotating the vector field if necessary, we can assume that $ \hat{y}^* < \infty $. Then, $ (0,\hat{y}^*) $ is a singular point of the $ x $-directional blown-up system $ \hat{X} $ and any of its asymptotic orbits correspond to asymptotic orbits of $ X $, with the exception $ \hat{x} = 0 $. Further, the assumptions that $ X $ has an isolated singularity and that the directional polynomial has only roots of multiplicity one implies that $ (0,\hat{y}^*) $ is a hyperbolic singular point, as argued in the proof of Proposition \ref{prop:DirectionImpliesOrbit}. But then, in order to have only one asymptotic orbit, $ (0,\hat{y}^*) $ must be a hyperbolic saddle and hence the eigenvalues must be of different sign. This results in the condition $ r^* < 0 $.
			
			Conversely, if the directional polynomial has precisely one root $ \hat{y}^* $, then by Proposition \ref{prop:DirectionImpliesOrbit} any asymptotic orbit corresponding to the singular point $ (0,\hat{y}^*) $ of system $ \hat{X} $ must also be an asymptotic orbit of $ X $ with characteristic direction $ \hat{y}^* $. Since the condition $ r^* < 0 $ implies this singular point is a hyperbolic saddle it is concluded that there is precisely one ingoing and one outgoing asymptotic orbit. Hence, by Theorem \ref{thm:C0Regularisable} we conclude that $ X $ is $ C^0 $-regularisable.
		\end{proof}
		
		Theorem \ref{thm:C0RegularisationHomX} removes a large class of vector fields that can not be block regularised. It provides as a corollary an easily computable necessary condition.
		
		\begin{corollary}\label{cor:C0regularisation}
			Let $ X\in\mathcal{V}^{hom} $. If $ X $ is $ C^0 $-regularisable then it is necessary that $ X $ is of homogeneity $ s $ even.
		\end{corollary}
		\begin{proof}
			A homogeneous vector field of degree $ s $ produces a directional polynomial of degree $ s+1 $. From Theorem \ref{thm:C0RegularisationHomX} it is necessary that there is precisely one real simple root, implying that the directional polynomial must be of odd degree. Hence, $ X $ must have homogeneity even.
		\end{proof}
	
		\begin{example}[Toy]
			From the computations in Example \ref{exam:desing} we know that the toy example $ X $ desingularises to a critical manifold that contains a unique hyperbolic singular point with eigenvalues $ \lambda_1 = (C_1 C_2^2)^{1/3}, \lambda_2 = -3(C_1 C_2^2)^1/3 $. Hence $ X $ contains one critical value $ r^* = \lambda_1/\lambda_2 = -1/3 $. As $ r^* < 0  $ we conclude from Theorem \ref{thm:C0RegularisationHomX} that $ X $ is $ C^0 $-regularisable.
		\end{example}
		
	\subsection{Stability of $ C^0 $ Regularity}\label{sec:stability}
		
		Let $ X_{hom} $ be the leading order homogeneous component of a smooth planar vector field $ X$, $X(0)=0 $. We want to know whether $ X_{hom} $ is $ C^0 $-regularisable if and only if $ X $ is $ C^0 $-regularisable, that is, whether it is enough to consider only $ X_{hom} $ in order to determine whether $ X $ is $ C^0 $-regularisable. As the argument for $ C^0 $-regularisation in Theorem \ref{thm:C0Regularisable} was purely topological, this question is answered if $ X $ is topologically equivalent to $ X_{hom} $. For $ X $ with leading order homogeneous component linear and hyperbolic, the well known Hartman-Grobman theorem confirms topological equivalence. Thankfully, this result has been extended by a theorem of Brunella and Miari \cite{Brunella1990}. The theorem is more general, dealing with quasi-homogeneous components. We restate the theorem in the special case of homogeneous components.
		
		\begin{thm}[Brunella-Miari]\label{thm:topologicalEquivalence}
			Suppose that the leading homogeneous component $ X_{hom} $ of $ X $ has an isolated singularity at $ 0 $. Then, in a neighbourhood of $ 0 $,
			\begin{enumerate}[i.]
				\item if $ X_{hom} $ has characteristic orbits or is a focus then $ X $ is topologically equivalent to $ X_{hom} $.
				\item if $ X_{hom} $ is a centre then $ X $ is a centre or focus.
			\end{enumerate}
		\end{thm}
		
		In Section \ref{sec:BlockRegularisation}, it was argued that any focus or centre is immediately not block regularisable. Hence, if $ X_{hom} $ is a focus or a centre and therefore not regularisable, then Theorem \ref{thm:topologicalEquivalence} shows $ X $ is also. Alternatively, if $ X_{hom} $ has characteristic orbits then Theorem \ref{thm:topologicalEquivalence} guarantees, by topological conjugacy, that $ X $ must have the same number of characteristic orbits. As characteristic orbits are also asymptotic orbits, we can show the following result.
		
		\begin{thm}\label{thm:C0stability}
			Let $ X_{hom} $ be the leading order homogeneous component of a planar vector field $ X $, $ X(0)=0 $. Suppose that $ X_{hom} $ has an isolated singularity at the origin. Then $ X $ is $ C^0 $-regularisable if and only if $ X_{hom} $ is $ C^0 $ regularisable.
		\end{thm}
		
		The preceding theorem applies to a generic class of vector fields. However, we have been concerned with vector fields $ X_{hom} \in \mathcal{V}^{hom} $. In what follows we wish to extend results about $ \mathcal{V}^{hom} $ to the following set of vector fields.
		\begin{definition}\label{def:V}
			Let $ X \in \mathcal{V} $ if it is a smooth planar vector field with leading order homogeneous component $ X_{hom} \in \mathcal{V}^{hom} $, i.e.{} without double roots, see Definition \ref{def:Vhom}. 
		\end{definition}
		Not only, from Theorem \ref{thm:C0stability} can it be concluded that vector fields $ X \in \mathcal{V} $ have their $ C^0 $-regularisation determined from their homogeneous component, but, Proposition \ref{prop:oneblowup} guarantees that they are the only set that requires precisely one blow-up to a desingularised vector field with only hyperbolic singularities on the critical manifold. This is crucial for ease of computation of the regularising map in what follows, but, the theory for the semi-hyperbolic case does exist and one could extend the the results in this paper to include all smooth vector fields covered in the hypothesis of Theorem \ref{thm:C0stability}.

%% file: RegularisingMap.tex
\section{$ C^k $ Regularisation}\label{sec:CkRegularisation}
	With $ C^0 $-regularisation firmly understood by topological properties of the vector field, we now begin the study of $ C^k $-regularity of the block map. We derive the formal asymptotic structure of the block map $ \pib $ for the class of vector fields $ \mathcal{V} $ defined in \ref{def:V}. It is shown that the block map can be explicitly computed as the composition of a smooth map conjugated by transition maps near hyperbolic saddles. These transition maps are well studied in, for example \cite{roussarie1995bifurcations}, and called \textit{Dulac maps}. Using the asymptotic series of the Dulac maps, necessary and sufficient conditions for $ C^k $-regularisation are proved.

	\subsection{Splitting the Problem}
		
		In Theorem \ref{thm:C0Regularisable} it was shown any vector field $ X \in \mathcal{V} $ is $ C^0 $-regularisable if and only if it has one ingoing and one outgoing asymptotic orbit. It is easily checked that $ X \in \mathcal{V} $ implies the union of the asymptotic orbits is smooth across the singularity. Assume that the vector field has been rotated so these asymptotic orbits are tangent to the $ x $-axis at $ 0 $. The rotation forces the characteristic direction $ \hat{y}^* = 0 $ and as such has no image in the $ y $ chart. Essentially, the blow-up and desingularisation produces two vector fields $ \hat{X} $ and $ \bar{X} $ that have a hyperbolic and regular point at $ 0 $, respectively. 
		
		Now, as there is only one asymptotic orbit, Lemma \ref{lem:Sections} implies that there is no need to study the block map on the entirety of an isolating block surrounding the singularity. Instead, it is sufficient to consider only two transversal sections $ \Sigma_0 \subset n^+$ and $ \Sigma_3 \subset n^- $ containing $ a^+ $ and $ a^- $ respectively.
		
\begin{figure}[ht]
			\centering
			\includegraphics[width=0.5\textwidth]{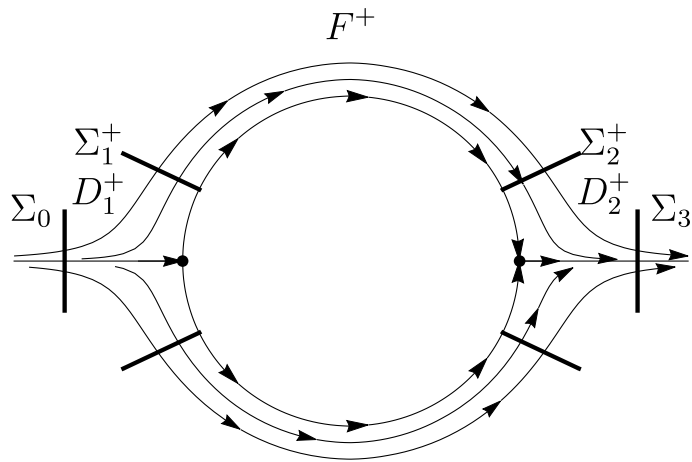}
			\caption{A schematic for the polar blow-up. The inner circle is the critical manifold $ r = 0 $.}\label{plot:polarblowup}
		\end{figure}

		The procedure is best conceptually understood through the polar blow-up, exampled in Figure \ref{plot:polarblowup}. To compute the block map $ \pib $ between the two transversal sections $ \Sigma_0 $ and $ \Sigma_3 $, we can introduce intermediate sections $ \Sigma_1 $ and $ \Sigma_2 $. The idea is to use the hyperbolicity and the regularity of the desingularised vector fields to explicitly find the form of the transitions between $ \Sigma_0 $ and the first intermediate section, between the two intermediate sections and then from the final intermediate section to $ \Sigma_3 $. Composing these transitions will yield the block map $ \pib $. 
		
\begin{figure}[ht]		
			\centering
			\input{DirectionalBlowUp.tikz}
			\caption{Intermediate sections and their desingularisations for the upper block map $ \pib_+ $.}\label{fig:DirectionalBlowUp}
		\end{figure}
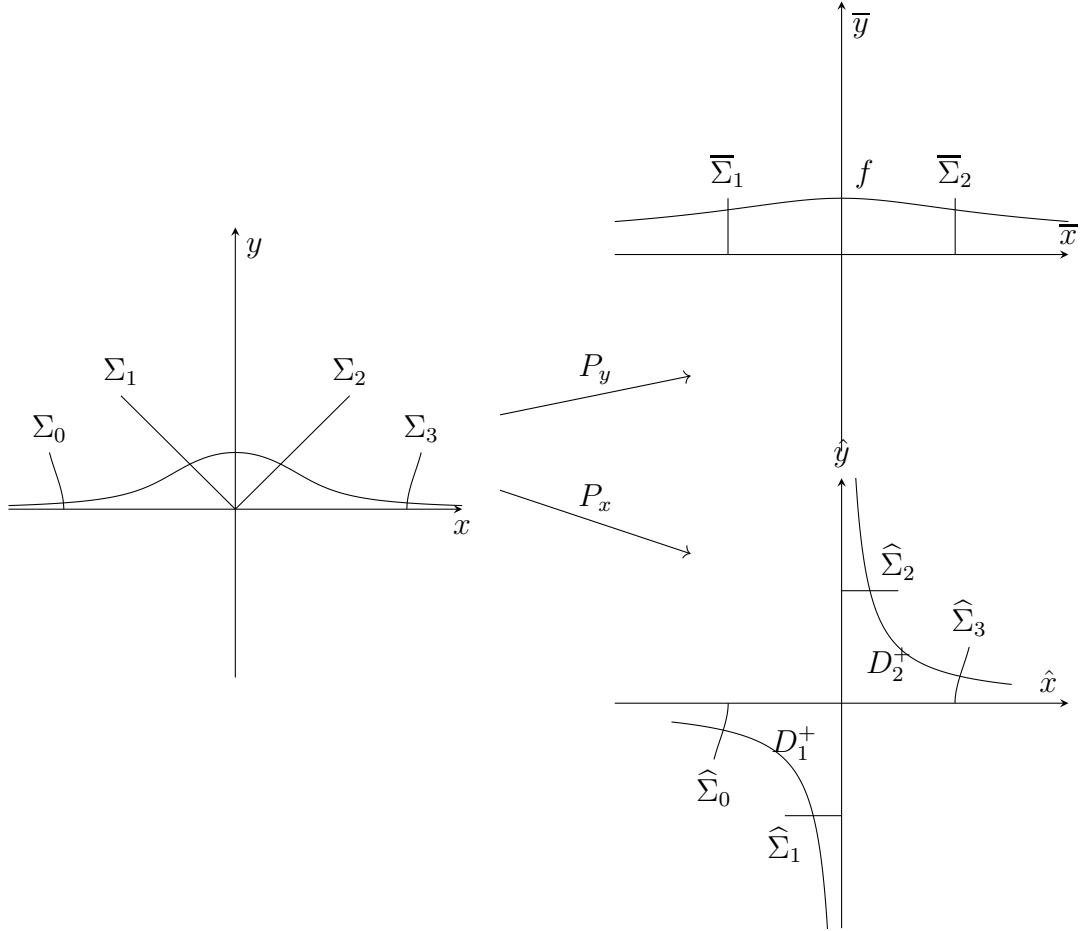

		Practically however, it is simpler to deal with the directional blow-up due to the corresponding vector fields being polynomial. Here $ \Sigma_1 $ and $ \Sigma_2 $ have images in $ \hat{X} $ and $ \bar{X} $, given by $ \hat{\Sigma}_i $ and $ \bar{\Sigma}_i $ for $ i = 1,2 $ respectively (see Figure \ref{fig:DirectionalBlowUp}).
		
		Note that the blow-up transforms $ P_x $ and $ P_y $ are discontinuous on the axes. This forces the separate study of the block map above and below the asymptotic orbit, denoted $ \pib_+ $ and $ \pib_- $. As a requirement for regularisation, the two maps must agree at $ 0 $ up to the expected order of differentiability. The upper half block map $ \pib_+ $ is displayed in Figure \ref{fig:DirectionalBlowUp}.
		
		More explicitly, if the hyperbolic transitions are denoted by $ D_1^\pm, D_2^\pm $ with $ \pm $ determining whether it is part of the upper or lower transition, and the regular transition is given by $ f:\overline{\Sigma}_1 \to \overline{\Sigma}_2 $, then the upper and lower block maps are given by
		\begin{align}
		\pib_+ &= D_2^+\circ F  \circ D_1^+, \label{eqn:P+}\\
		\pib_- &= (D_1^-)^{-1} \circ F^{-1}\circ (D_2^-)^{-1}\label{eqn:P-}
		\end{align}
		where $ F = P_2^{-1} \circ f \circ P_1 $ and $ \displaystyle \left. P_i = P_y\circ P_x^{-1}\right|_{\hat{\Sigma}_i}$, see Figure \ref{fig:DirectionalBlowUp}.
		
		The next step in the derivation is a discussion of the form of the hyperbolic transitions $ D_i^\pm $. 
		
	\subsection{Transition Near a Hyperbolic Saddle}
		Most of the work in this subsection is contained in \cite{roussarie1995bifurcations}. Consider a hyperbolic singular point with eigenvalues $ \lambda_2 < 0 < \lambda_1 $. Let $ \displaystyle r = -\frac{\lambda_1}{\lambda_2} $ be the \textit{ratio of hyperbolicity}. Up to translation and rescaling, the vector field $ X $ takes the form,
		\begin{equation}\label{eqn:nonNormal}
			X = \left( x + \dots, -r y + \dots \right).
		\end{equation}
		The following theorem is a version of the Dulac-Poincar\'{e} normal form.
		\begin{thm}[\protect{\cite{roussarie1995bifurcations}}]\label{thm:NormalForm}
			If $ r = p/q \in \Q $ with $ p,q $ co-prime, there exists a $ C^\infty $ change of variables $ \phi $ and a rescaling $ K $ transforming  $ X $ to the normal form,
			\begin{equation}\label{eqn:normal}
				X_N = \left(x, -r y + y\frac{1}{q} \sum_{i = 1}^{\infty} \alpha_{i+1} (x^p y^q)^i\right)	
			\end{equation}
			If $ r \nin \Q $ then $ \alpha_{i+1} \equiv 0 $ for $ i \geq 1 $. In other words, $ X $ is (orbitally) $ C^\infty $ conjugate to the \textit{normal form} $ X_N $.
		\end{thm}
		A proof can be found in \cite{Belitskii2002} along with the iterative procedure of computing $ \phi $. 
%		We remark that the proof involves several steps. The first is showing that there is a formal transformation $ \hat{\phi} $ bringing the system to normal form. Then, using a modification of the Borel lemma, it is shown that there is a $ C^\infty $ transformation $ \bar{\phi} $ asymptotic to $ \hat{\phi} $ bringing $ X $ to $ X_N + T $, where $ T $ is a flat residual (i.e. has vanishing derivatives of all order at $ 0 $). Then, the Sternberg-Chen theorem \cite{Belitskii2002} guarantees, if $ X $ is hyperbolic, there is another $ C^\infty $ transformation killing the residual. The composition of this transformation with $ \bar{\phi} $ produces the desired $ \phi $.

		Now, consider the two transversal sections $ \sigma = [0,1) \times \{1\} $ and $ \tau = \{1\}\times [0,1) $ near the hyperbolic point (see Figure \ref{fig:Dulac}). We wish to study the structure of the transition map $ D:\sigma \to \tau $, denoted the \textit{Dulac map}, for the normal form $ X_N $. In particular we seek the asymptotic expansions near $ x = 0 $.
		
		\begin{figure}[H]
			\centering
			\begin{tikzpicture}[scale = 4]
			\draw[->- = 0.5] (0,0) to (1.1,0);
			\draw[->- = 0.5] (0,1.1) to (0,0);
			\node[left] at (0,1) {$ y $};
			\node[below] at (1,0) {$ x $};
			\draw[domain=0.1:1,->- = 0.5] plot (\x,0.1/\x);
			\draw[domain=0.2:1,->-=0.5] plot (\x,0.2/\x);
			\draw[domain=0.3:1,->-=0.5] plot (\x,0.3/\x);
			\draw (0,1) -- (0.5,1) node[above,at end] {$ \sigma $};
			\draw (1,0) -- (1,0.5) node[right,at end] {$ \tau $};
			\end{tikzpicture}
			\caption{}\label{fig:Dulac}
		\end{figure}
		
		Firstly, from Figure \ref{fig:Dulac} it is clear that $ D(x) $ admits a continuous extension to $ D(0) = 0 $. Using the flow-box theorem, it can be reasoned that $ D $ is smooth in $ x $ for $ x \neq 0 $ if the vector field $ X $ is smooth. Now, if $ r \nin \Q $ then one can integrate $ X_N $ to see $ D $ is simply $ D(x) = x^r $. However, if $ r = p/q $ the asymptotic expansion of $ D $ becomes significantly more complicated.
		
		The following proposition collects numerous results from \cite{roussarie1995bifurcations} and describes the asymptotic expansion near $ x = 0 $ of $ D $. A small part of the proof is included as it is constructive and useful for computations in Section \ref{sec:RegularisingMap}.
		\begin{proposition}[\protect{\cite{roussarie1995bifurcations}}]\label{prop:normalform}
			If $ r = p/q \in \Q $ there exists a sequence of polynomials $ P_k $ of degree $ k-1 $ such that the formal series
			\[ \hat{D}(x) = x^r + \sum\limits_{k=2}^{\infty} x^{k p+r} P_k(\ln x), \]
			is asymptotic to $ D $ in the sense that, for any $ s $,
			\[ |D(x) - \sum_{k=1}^{s} x^{kp+r}P_k(\ln x) | = O(x^{s p +r}). \]
			If $ r \nin \Q $ then $ \hat{D}(x) = x^r $.
		\end{proposition}
		\begin{proof}
			The formula is a direct result from the normal form $ X_N $. First, make the singular change of variables $ u = x^p y^q,\ x = x $ to obtain the decoupled equations,
			\begin{equation}\label{eqn:normalform}
				\begin{aligned}
					\dot{x}	&= x \\
					\dot{u}	&= \sum_{k = 1}^{\infty} \alpha_{k+1} u^{k + 1} = P(u).
				\end{aligned}
			\end{equation}
			The first equation gives $ x = x_0 e^t $ providing the transition time from $ \sigma $ to $ \tau $ as $ T_x = - \ln x $. In the second equation $ P(u) $ is analytic in some neighbourhood of $ u = 0 $. For each $ t $ we can expand the series in $ u_0 $ to obtain the solution,
			\begin{equation}
			u(t,u_0) = \sum_{k=1}^{\infty} g_k (t) u_0^k.
			\end{equation}
			The $ g_k $ can be found as the solutions to higher order variational equations with $ g_1(0) = 1 $ and $ g_k(0) = 0 $ for $ k \geq 2 $. In particular, $ g_1(t) = 1 $.
			
			Hence, using the fact that at $ t = 0 $, $ u = x^p $ and at $ t = T_x $, $ u = y^q $ then 
			\begin{align*}
			y^q		&= \sum_{k=1}^{\infty} g_k (T_x) u_0^k	\\
			\hat{D}(x)^q	&= x^p + \sum_{k=2}^{\infty} g_k (-\ln x) x^{k p}.
			\end{align*}	
			The fact that $ g_k $ is polynomial of degree $ k-1 $ is proved in \cite[p.~99]{roussarie1995bifurcations}.	
		\end{proof}
		The proposition motivates the following definition of quasi-regularity.
		\begin{definition}
			A germ of a map $ f $ at $ 0\in \R^+ $ is called \textit{quasi-regular} if
			\begin{enumerate}[i.]
				\item $ f $ has a representative on $ [0,a) $ that is $ C^\infty $ on $ (0,a) $.
				\item $ f $ is asymptotic to the Dulac series,
				\[ \hat{f}(x) = \sum_{k=1}^{\infty} x^{\lambda_k}P_k(\ln x), \]
				with $ (\lambda_1,\lambda_2,\dots) $ a positive, strictly increasing sequence tending to infinity and  $ P_k $ a sequence of polynomials. We say $ f $ is a \textit{quasi-regular homeomorphism} if $ f $ is quasi-regular and $ P_1(x) \equiv c \in \R^+ $.  
			\end{enumerate}
		\end{definition}
		Note that these maps are also known as \textit{almost regular} in the literature \cite{Ilyashenko2008}. 
		
	\subsection{Quasi-Regularity of the Block Map}\label{sec:RegularisingMap}
		
		It can now be argued that the block map $ \pib $ is quasi-regular and hence, any $ C^0 $-regularisable $ X\in\mathcal{V} $ generically has only finite regularity. We remark that the set of all quasi-regular homeomorphisms $ \mathcal{D} $ is a group under composition. Further, the group of diffeomorphisms fixing $ 0 $ is a subgroup of $ \mathcal{D} $. 
		
		Now, consider any of the $ D_i^{\pm} $, for instance $ D_1^+ $. A method of computation for $ D_i^{\pm} $ will be constructed in section \ref{ssec:iterativemethod} but a quick description is as follows. Construct the $ C^\infty $ normalising transformation bringing $ \hat{X} $ into the normal form $ X_N $. The transform exists due to Theorem \ref{thm:NormalForm}. Flow the image $ \phi\circ\hat{\Sigma}_0 $ to $ \tau = \set{1}\times[0,1) $. Then the asymptotic series of the Dulac map $ D:\tau \to \sigma $, with $ \sigma = [0,1)\times\set{1} $, can then be computed according to Proposition \ref{prop:normalform}. Lastly, flow $ \sigma $ to $ \phi\circ\hat{\Sigma}_1 $.
		
		Therefore, any of the $ D_i^{\pm} $ are obtained by the composition of smooth maps and a Dulac map $ D $. From the group structure of $ \mathcal{D} $, the following proposition is immediate and provides exactly what we have searched for; the finite differentiability of the block map.
		
		\begin{proposition}
			The upper and lower block maps $ \pib_{\pm} $ are quasi-regular.
		\end{proposition}
		
		It can now be seen that there are two mechanisms withholding $ \pib $ from being $ C^\infty $. The first is that both $ \pib^\pm $ are in general only quasi-regular and thus may contain terms of type $ x^{\lambda}(\ln x)^k $, forcing at best $ \pib \in C^{\lambda-1} $. The second is the potential for the two quasi-regular maps $ \pib^\pm $ to disagree for some terms in the series. The rest of this section will explore these possibilities in more detail.
		
		Let $ X_N $ be the normal form $ C^\infty $ equivalent to $ \hat{X} $ and $\phi $ the normalising transformation. In what follows we make a particular choice of transversal sections $ \Sigma_0, \Sigma_3 $ and the intermediate sections $ \Sigma_1, \Sigma_2 $ in Figure \ref{fig:DirectionalBlowUp}. This choice is justified by noticing that away from the singular point $ 0 $ of $ X $, the flow-box theorem guarantees the transition between two smooth transverse sections in a close enough neighbourhood of each other is a diffeomorphism. Hence, any results about finite differentiability will hold for any nearby sections.
		
		The choice of sections are as follows. Let 
		\[ \tau_0 = \set{1}\times [0,1),\ \sigma_1 = [0,1) \times\set{1},\ \sigma_2 = (-1,0]\times\set{-1},\ \tau_3 = \set{1}\times [0,1), \]
		defined in the normal form system $ X_N $. Then, choose $ \Sigma_0, \Sigma_1, \Sigma_2, \Sigma_3 $ to be the pre-images under the map $ P_x \circ \phi $ of $ \tau_0, \sigma_1, \sigma_2, \tau_3 $ respectively. As described in Section \ref{sec:RegularisingMap}, let $ \hat{\Sigma}_i $ and $ \bar{\Sigma}_i $ be the images of $ \Sigma_i $ under $ P_y $ and $ P_x $ respectively (see Figure \ref{fig:DirectionalBlowUp}).
		
		The hyperbolic transitions $ D_i^\pm $ can be written as the composition $ \phi^{-1}\circ d_i^\pm $ where $ d_i $ are the hyperbolic transition maps in the normal form coordinates. Let $ r = -r^* $ be the negative of the critical value of $ X $ whose unique existence is guaranteed from Theorem \ref{thm:C0RegularisationHomX}. Then, as determined in Proposition \ref{prop:normalform}, $ d_i $ are Dulac maps of the form
		\begin{align}
			d_1^\pm(x)		&= x^{r} + O(x^{2p+r}\ln x),	\\
			d_2^\pm(x)		&= x^{1/r} + O(x^{2q+1/r}\ln x),
		\end{align}
		if $ r=p/q\in\Q $, $ p,q $ coprime. If $ r \nin \Q $ then $ d_i^{\pm} $ is equivalent to the leading order.
		
		The equations for the upper and lower block maps (\eqref{eqn:P+}, \eqref{eqn:P-}) are,
		\begin{align}
			\pib_+ &= d_2^+\circ F \circ d_1^+,\\
			\pib_- &= (d_1^-)^{-1} \circ F^{-1}\circ (d_2^-)^{-1},
		\end{align}
		where $ F = P_2^{-1} \circ f \circ P_1  $ and $ P_1: \sigma_1 \to \bar{\Sigma}_1 $, $ P_2: \sigma_2 \to \bar{\Sigma}_2 $, are the two images of $ \sigma_i $ under the smooth transformations $ P_y\circ P_x^{-1}\circ\phi^{-1} $. As $ F $ is a smooth function, it can be assumed that the asymptotic expansion about $ 0 $ is a diffeomorphism of the form
		\[ F(x) = F^\prime(0) x + \frac{1}{2!}F^{\prime\prime}(0) x^2  \dots,\quad F^\prime(0) > 0. \]
		
		Firstly, it can be seen that the composition does not change the order of the leading order term of $ \pib^\pm $; a direct computation yields,
		\begin{equation}\label{eqn:firstvalPib}
			\begin{aligned}
				\pib_+(x) &=  \left(F^\prime(0)\right)^{1/r} x + \dots, \\
				\pib_-(x) &=  \left((F^{-1})^\prime(0)\right)^{1/r} x + \dots.
			\end{aligned}
		\end{equation}
		where $ +\dots $ represents terms that do not effect the first one-sided derivatives at $ 0 $. If it is required that $ \pib $ is a diffeomorphism, i.e.{} at least $ C^1 $,
		then we must insist that $ \pib_+ $ and $ \pib_- $ agree at 0 up to the first derivative. This gives rise to the following lemma. 
		
		\begin{thm} \label{thm:NearId}
			 $ X\in\mathcal{V} $ is $ C^1 $-regularisable if and only if $ F^\prime(0) = 1 $. In this case, the block map $ \pib $ is in fact a diffeomorphism. Else, $ X $ is $ C^{0,1} $-regularisable.
		\end{thm}
		\begin{proof}
			The proof is a simple computation. In order for $ \pib(x) $ to be a diffeomorphism about $ 0 $ it is necessary that the following string of implications holds.
			\begin{align*}
				\pib_+^\prime(0) = \pib_-^\prime(0) \iff F^\prime(0) &= (F^{-1}(0))^\prime \\
				\iff F^\prime(0) &= (F^\prime(0))^{-1}	\\
				\iff F^\prime(0) &=  1. 
			\end{align*}
		\end{proof} 
		
		The following definitions will enable us to understand the next highest order term in the asymptotic series of $ \pib^\pm $.
		\begin{definition}
			If $ r \in \Q $ then we say the saddle point is of finite order if there exists a smallest integer $ m \geq 2 $ such that the coefficient of $ (x^p y^q)^{m-1} $ in the normal form \eqref{eqn:normal}, $ \alpha_m $, is non-zero. That is, $ P(u) $, defined in \eqref{eqn:normalform}, takes the form,
			\begin{equation}
				P(u) = \alpha_m u^m + o(u^m),\quad \alpha_m \neq 0.
			\end{equation}
			The integer $ m $ is called the order of the resonance.
		\end{definition}
		
		With this definition, a corollary of Proposition \ref{prop:normalform} follows.
		\begin{corollary}\label{cor:finiteResonanceDulacMap}
			Suppose that a saddle point has a finite order of resonance $ m $. Then, the Dulac map and it's inverse are given by,
			\begin{equation}
			\begin{aligned}
				D(x) &= x^r \left(1 - \frac{1}{q} \alpha_m x^{(m-1)p}\ln x + \dots \right)\\
				D^{-1}(y) &= y^{1/r} \left(1 + \frac{q}{p^2} \alpha_m  y^{(m-1)q}\ln y + \dots \right).	
			\end{aligned}
			\end{equation}
		\end{corollary}
		
		\begin{definition}
			The \textit{order of transition} of a diffeomorphic transition map $ F $ is the first integer $ k >1 $ such that $ F^{(k)}(0)\neq 0 $. That is, the map is of the form,
			\[ F(x) = F^{(1)}(0) x + \frac{1}{k!}F^{(k)}(0) x^k + \dots. \]
		\end{definition}
		
		With these definitions the first higher order term of $ \pib^\pm $ can now be expressed. Suppose $ X_N $ has resonance of order $ m $ and the order of the transition map $ F $ is $ k $. Further suppose that $ X $ is $ C^1 $-regularisable, else no investigation of higher order terms would be required to determine the regularity of $ X $. Then,
		\begin{align*}
			d_1^\pm(x) 	&= x^r \left(1 - \frac{1}{q} \alpha^\pm_{1,m} x^{(m-1)p}\ln x + \dots \right) \\
			d_2^\pm(x) 	&= x^{1/r} \left(1 - \frac{1}{p} \alpha^\pm_{2,m}   x^{(m-1)q}\ln x + \dots \right) \\
			F(x)		&= x + \frac{1}{k!}F^{(k)}(0) x^k + \dots. 
		\end{align*}
		Composing the maps achieves the following lemma.
		\begin{lemma}\label{lem:leadingorder}
			Let $ X $ be $ C^1 $-regularisable. Then to leading higher order,
			\begin{equation}
			\begin{aligned}
				\pib^\pm(x) &=  x\left( 1+\psi^\pm(x) + \dots\right),\\ \psi^\pm(x) &= \begin{cases}
					\displaystyle\pm\frac{1}{k! r} F^{(k)} x^{(k-1)r}, & (k-1) < (m-1)q \\[2ex]
					-\left( \frac{1}{p} \alpha^\pm_{1,m} + \frac{1}{q} \alpha^\pm_{2,m} \right)x^{(m-1)p}\ln x, & (k-1) \geq (m-1)q
				\end{cases}
			\end{aligned}
			\end{equation}
			with $ F^{(k)} := F^{(k)}(0) $.
		\end{lemma}
		It is now clear how the loss of regularity arises. The following mechanisms cause a loss of regularity.
		\begin{enumerate}
			\item $ F^\prime(0) \neq 1 $ producing a discrepancy in the upper and lower block maps. Essentially, orbits one one side are 'pulled in` whilst orbits on the other side are 'pushed out`. In such a case $ X $ is $ C^{0,1} $-regularisable.
			\item The transition $ F $ produces a term $ F^{(k)} \neq 0 $ before any resonance terms from the saddle. This corresponds to the first case of $ \psi $ and can be seen to create a term in the series of order $ x^{(1+(k-1)r)} $. Hence, $ X $ is $ C^{1+(k-1)r} $-regularisable, provided $ r\nin\N $. If $ r \in \N $, then the difference in sign between $ \pib^\pm $ forces $ X $ to be $ C^{(k-1)r,1} $-regularisable.
			\item A resonance term is produced affecting the series prior to any discrepancies in the transition $ F $. This is the second case of $ \psi $. Provided $ \alpha^\pm:=\left( \frac{1}{p} \alpha^\pm_{1,m} + \frac{1}{q} \alpha^\pm_{2,m} \right) \neq 0 $, the block map $ \pib $ is $ C^{(m-1)p} $ and $ \pib^{((m-1)p)} $ is Log-Lipschitz continuous. This implies $ X $ is $ C^{(m-1)p,a} $-regularisable, for some 
			$ 0 \leq a < 1 $. If, for example, $ d_2^\pm = (d_1^\pm)^{-1} $, then it may occur that $ \alpha^\pm = 0 $ and higher order terms will need to be computed to determine the regularity.
		\end{enumerate}

		Each of the possible mechanisms are exemplified in Section \ref{sec:examples} and phase portraits of these examples are given in Figure \ref{fig:examples}. The problem of finite regularity has been reduced to a computation of the order of transition $ k $ and the order of resonance $ m $. In the next section some of the difficulties in computing these quantities are overcome. 	
	
	\subsection{A Method for Computing the Block Map}\label{ssec:iterativemethod}
		In the previous section the mechanisms for finite regularity were determined. To see if finite regularity occurs for some $ C^0$-regularisable $ X \in \mathcal{V} $, the Dulac maps $ d_i^\pm $ and the smooth transition map $ F $ must be computed. Whilst there is much literature on computing the hyperbolic transition $ d_i $, for instance \cite{roussarie1995bifurcations}, an issue arises when trying to compute the transition between the saddles, $ F $. The transition is between $ \sigma_1 $ and $ \sigma_2 $ in the normal form $ X_N $. However, in order to calculate the map, local coordinates in $ X_N $ need to be transformed to non-local coordinates in $ \bar{X} $. Here lies the problem; at any iteration in the computation of $ \phi $, $ \phi $ is only known up to some truncated order. So the image of $ \sigma_i $ in $ \bar{X} $ will only be known to some truncated order and hence, when solving the variational equations to compute $ F $, some terms in the varied coefficients will not have yet been computed. 		
		
		The key to resolving the problem is to observe that the block map $ \pib $ is independent of the choice of intermediate sections $ \Sigma_1 $ and $ \Sigma_2 $. Hence, there is freedom in the choice of these sections. In particular, the sections can be chosen in $ X_N $ as \[ \sigma_1^\epsilon = [0,1)\times\set{\epsilon},\qquad \sigma_2^\epsilon = (-1,0]\times\set{-\epsilon}, \]
		and $ \Sigma_i^\epsilon,\hat{\Sigma}_i^\epsilon,\bar{\Sigma}_i^\epsilon $ the image of $ \sigma_i^\epsilon $ in $ X, \hat{X},\bar{X} $ respectively. As the final series does not depend on this choice, then it must be that the series does not depend on $ \epsilon $, and inevitably, we can send $ \epsilon \to 0 $.
		
		The steps in the method needed to compute $ F $ up to order $ k $ are given below. The computation of the first order term $ F^{(1)} $ is done as an example. It is assumed that the characteristic direction of $ X $ is $ 0 $.
		\begin{enumerate}
			\item Truncate $ X $ at order $ k $ to produce the vector field \[ X_k = \left(P_s(x,y) + \dots + P_{k+s}(x,y),\ Q_s(x,y) + \dots + Q_{k+s}(x,y)\right), \]
			where $ P_i,Q_i $ contain the homogeneous components of order $ i $.
			\item Compute the blow-ups $ \hat{X}_k,\bar{X}_k $:
			\begin{align*}
				\hat{X}_k &= \left(\hat{x} P_s(1,\hat{y}) + \dots + \hat{x}^{k+1} P_{s+k}(1,\hat{y}),\ p_s(\hat{y}) + \dots +  \hat{x}^k p_{s+k}(\hat{y})\right) \\
				\bar{X}_k &= \left(q_s(\bar{x}) + \dots + \bar{y}^k q_{s+k}(\bar{x}),\ \bar{y} Q_s(1,\bar{x}) + \dots +  \bar{y}^{k+1} Q_{s+k}(1,\bar{x})\right),   
			\end{align*}			
			with $ p_i(\hat{y}) = Q_i(1,\hat{y}) - \hat{y}P_i(1,\hat{y}) $ and $ q_i(\bar{x}) = P(\bar{x},1) - \bar{x}Q(\bar{x},1) $ the higher order directional polynomials in the $ x $ and $ y $ charts respectively. 
			\item Compute the normal form $ X_N $ and normalising transformation $ \phi:X\to X_N $ to order $ k+1 $.
			
			\enumeratext{For the first order term in $ F $ we have that $ k=0 $, $ \phi =(x,y)+||(x,y)||^2 $, $ X_N = (x,-r y)+||x,y||^2 $ and $ r $ is the characteristic value of $ X $.}
			
			\item Using $ P_y\circ P_x^{-1}\circ\phi^{-1}:X_N\to\bar{X} $ compute to order $ k+1 $ the images of $ \sigma_1^\epsilon,\sigma_2^\epsilon $ in $ \bar{X}_k $, that is $ \bar{\Sigma}_1^\epsilon,\bar{\Sigma}_2^\epsilon $. 
			
			\enumeratext{For the leading order term this produces
				\begin{align*}
					\bar{\Sigma}_1^\epsilon &= \set{((1/\epsilon + O(\epsilon)) +O(x_1),(\epsilon+O(\epsilon^2)) x_1 + O(x_1^2))| x_1\in[0,1)} \\
					\bar{\Sigma}_2^\epsilon	&=  \set{((-1/\epsilon + O(\epsilon)) +O(x_2),(\epsilon+O(\epsilon^2)) x_2 + O(x_2^2))| x_2\in[0,1)}
				\end{align*}
			}
			\item Use the $ (k+1)^{th} $ order variational equations to find $ f_\epsilon:\bar{\Sigma}_1^\epsilon\to \bar{\Sigma}_2^\epsilon $ and then use the preimage under $ P_y\circ P_x^{-1}\circ\phi^{-1} $ to obtain $ F_\epsilon $ to order $ k+1 $.
			
			\enumeratext{By definition $ F(x_1) = x_2(x_1) $. Using the formula from \cite[p.~68]{andronov1974qualitative}, the transition to first order is given by
			\[ F(x_1) = (1+O(\epsilon))\exp\left(\int_{-1/\epsilon+O(1)}^{1/\epsilon+O(1)} \frac{Q_s(x,1)}{q_s(x)} dx  \right) x_1 + O(x_1^2). \]   }
		
			\item Take the limit $ \epsilon\to 0 $ to obtain $ F $ to order $ k $.
		\end{enumerate}
		
		The limit in step 6 can be seen as the same calculation involved in computing the Cauchy Principal Value of the integral in step 5. Hence, we have the following nice corollary of Theorem \ref{thm:NearId}.
		\begin{corollary}\label{cor:Cauchy}
			Let $ X \in\mathcal{V} $ be $ C^0 $-regularisable and normalised so that the asymptotic orbit is tangent to the $ x $-axis, that is, $ \hat{y}^*=0 $. Then $ X $ is $ C^1 $-regularisable if and only if the Cauchy Principal Value of 
			\begin{equation}\label{eqn:C1reg}
			\int_{-\infty}^{\infty} \frac{ Q_s(x,1)}{q_s(x)} dx = 0.
			\end{equation} 
			Else, $ X $ is $ C^{0,1} $-regularisable.
		\end{corollary}	
		As $ \displaystyle \frac{Q_s(1,x)}{q_s(x)} $ is a rational function, this is equivalent to the requirement that the sum of its residues in the upper half complex plane be zero. Note that if a singularity $ x^* \in \R\cup\set{\infty} $ then we must take half the residue in the sum as it lies on the contour. The beauty of this corollary is that the $ C^1 $ regularity of $ X $ can be computed purely from the leading order terms in $ X $, with no need for blow-up, normal forms or variational equations.
		
		\begin{example}[Toy]\label{exam:toyC1}
			Recall the toy example
			\[ X = \left( C_1 z_2^2,\ C_2 z_1^2 \right). \]
			In order to normalise $ X $ to use Corollary \ref{cor:Cauchy} make the transformation $ (z_1,z_2)\mapsto(C_1^{1/3} z_1, C_2^{1/3} z_2) $ and rescale by $ (C_1 C_2)^{-2/3} $ to get the vector field,
			\[ X^\prime = \left( z_2^2,\ z_1^2 \right). \]
			Now we align the asymptotic orbit with the $ z_1 $-axis by rotating by $ \pi/4 $ and achieving, up to rescaling,
			\[ X^{\prime\prime} = \left( z_1^2 + z_2^2,\ -2 z_1 z_2 \right). \]
			Then $ Q_2(z_1,1) = -2 z_1 $ and $ q_2(z_1) = P_2(z_1,1)-z_1 Q_2(z_1,1) = 3 z_1^2 +1. $ From Corollary \ref{cor:Cauchy} $ X $ is $ C^1 $ if and only if the residues of 
			\[ g(z):= \frac{ Q_2(z_1,1)}{q_2(z_1)} = -\frac{2 z_1}{3 z_1^2+1}, \] sum to $ 0 $. We have singularities in $ g(z) $ at $ \pm 1/\sqrt{3} i,\infty $. Taking only the singularities in the upper half complex plane, the Cauchy principal value is \[ \int_{-\infty}^\infty g(z) dz = \Res_{z=\frac{1}{\sqrt{3} i}}\left(g(z)\right) + \frac{1}{2} \Res_{z=\infty}(g(z)) = -\frac{1}{3}+\frac{1}{3} = 0. \] We conclude that the toy example $ X $ has at least $ C^1 $ regularity.
		\end{example}
		

%% file: DirectionalBlowUp.tikz
\begin{tikzpicture}
				\begin{axis}[name = plot1,
					scale only axis,
					width = 1/3*\textwidth,
					height =1/3*\textwidth ,
					axis x line=middle,
					axis y line=middle,
					x label style = {anchor=north},
					axis equal,
					xlabel = {$x$},
					ylabel = {$y$},
					restrict y to domain = -4:2,
					restrict x to domain = -2:2]
					\addplot [domain = 0:0.5, samples = 300]
					({sqrt(0.125/x -x^2)}, {x});
					\addplot [domain = 0:0.5, samples = 300]
					({-sqrt(0.125/x -x^2)}, {x}); 
					\addplot [domain = 0:0.5, samples = 300]
					(1.5+x^2-x^3, {x}) node[above,pos=1]{$ \Sigma_3 $};
					\addplot [domain = 0:0.5, samples = 300]
					(-1.5-x^2+x^3, {x}) node[above,pos=1]{$ \Sigma_0 $};
					\addplot [domain = 0:1, samples = 300]
					({x}, {x}) node[above,pos=1]{$ \Sigma_2 $};
					\addplot [domain = 0:-1, samples = 300]
					({x}, {-x})node[above,pos=1]{$ \Sigma_1 $};
				\end{axis}
				
				\begin{axis}[name = plot2,at={($(plot1.east)+(2cm,3cm)$)},anchor=outer west,axis x line=middle,
					scale only axis,
					width = 1/3*\textwidth,
					height =1/3*\textwidth ,
					axis y line=middle,
					x label style = {anchor = south},
					axis equal,
					xlabel = {$\overline{x}$},
					ylabel = {$\overline{y}$},
					restrict y to domain = -2:2,
					restrict x to domain = -2:2]
					\addplot [domain = 0:0.5, samples = 300]
					({sqrt(0.125/x^3 -1)}, {x}) node[above,pos=0.9]{$ f $};
					\addplot [domain = 0:0.5, samples = 300]
					({-sqrt(0.125/x^3 -1)}, {x});
					\addplot [domain = 0:0.5, samples = 300]
					(1, {x}) node[above,pos=1]{$ \overline{\Sigma}_2 $};
					\addplot [domain = 0:0.5, samples = 300]
					(-1, {x}) node[above,pos=1]{$ \overline{\Sigma}_1 $};  
				\end{axis}
				
				\begin{axis}[name = plot3,at={($(plot1.east)+(2cm,-3cm)$)},anchor=outer west,axis x line=middle,
					scale only axis,
					width = 1/3*\textwidth,
					height =1/3*\textwidth ,
					axis y line=middle,
					y label style = {anchor=south},
					axis equal,
					xlabel = {$\hat{x}$},
					ylabel = {$\hat{y}$},
					ymin = -2,ymax=2,
					xmin = -2, xmax = 2]
					\addplot [domain = 0.09:1.5, samples = 300]
					({x}, {0.25/x}) node[below,pos=0.65]{$ D_2^+ $}; 
					\addplot [domain = -0.09:-1.5, samples = 300]
					({x}, {0.25/x}) node[above,pos=0.65]{$ D_1^+ $}; 
					\addplot [domain = 0:0.5, samples = 300]
					(1+x^2-x^3, {x}) node[above,pos=1]{$ \widehat{\Sigma}_3 $};
					\addplot [domain = -0.5:0, samples = 300]
					(-1-x^2-x^3, {x}) node[below,pos=0]{$ \widehat{\Sigma}_0 $};
					\addplot [domain = 0:0.5, samples = 300]
					(x, 1) node[above,pos=1]{$ \widehat{\Sigma}_2 $};
					\addplot [domain = -0.5:0, samples = 300]
					(x, -1) node[below,pos=0]{$ \widehat{\Sigma}_1 $};
				\end{axis}
				
				\draw[->] ($(plot1.east) + (0.5cm,0.5cm)$) -- ($ (plot2.south west) +(1cm,1cm)  $) node[midway,above]{$ P_y $};
				\draw[->] ($(plot1.east) + (0.5cm,-0.5cm)$) -- ($ (plot3.north west) + (1cm,-1cm)$) node[midway,above]{$ P_x $};

			\end{tikzpicture}

%% file: Normalising.tex
\section{Examples of Finite Regularity}\label{sec:examples}
	With the theory and methods outlined, some examples of the mechanisms producing finite regularisation are produced in this section. In particular, a complete classification of the regularity of homogeneous quadratic vector fields in the generic class $ {V}^{hom} $ is given. Firstly, a canonical form for homogeneous vector fields is derived. The quadratic vector fields are then categorised using this canonical form and the regularity of each determined. Lastly, higher order perturbations of the quadratic vector fields are considered, showing higher order finite regularity.
	
	\subsection{Canonical Form for Homogeneous Polynomial Vector Fields}\label{sec:Normalisation}
		
		The aim of this section is to propose a canonical form for homogeneous polynomial vector fields in the plane. The proposed canonical form is justified geometrically and agrees with the Jordan canonical form when the vector field is of homogeneity 1. 
		
		Consider the set of planar vector fields $ X = (P(x,y), Q(x,y)) $ with $ P,Q $ homogeneous polynomials of degree $ n $. We are interested in the set of equivalence classes under linear transformations $ \mathrm{GL}(2,\R) $ of the coordinates $ (x,y) $. As evidenced in earlier sections, the directional polynomial $ p(x,y) = x Q(x,y) - y P(x,y) $ is a fundamental object of study for $ X $. The directional polynomial can be characterised by its roots $ \omega = [x^*,y^*] $ over real projective space $ \PP $, which have been denoted the characteristic directions. Our choice of a canonical form will be justified by considering the induced action of $ \mathrm{GL}(2,\R) $ on the characteristic directions.
		
		In the case of a hyperbolic saddle, the unstable and stable manifolds in canonical form are the $ x $ and $ y $ axis, thus setting the characteristic directions to $ \omega = [0,1],[1,0] $. We seek an extension of this process by suggesting a canonical form of the directional polynomial for arbitrary degree. Geometrically, a choice of canonical form fixes the characteristic directions and thus the asymptotic direction of the characteristic orbits. The problem of selecting a canonical form is reduced to the study of the equivalence classes of characteristic directions under linear transformations of the original coordinates.
		
		For what follows, assume the directional polynomial $ p(x,y) \neq 0 $ and consider $ \PP $ in the chart $ \omega = [1,v]=[1,y/x] $. Then, $ p(v) := p(1,v) $ has precisely $ n+1 $ roots in $ \hat{\C} :=  \C\cup \{\infty\} $, which is denoted by the set $ \hat{\omega} = [1,\hat{v}] \in \PP(\hat{\C})^{n+1} $. The following proposition describes how $ \hat{\omega} $ is transformed when $ x,y $ undergo linear transformations.
		
		\begin{proposition}
			The action of $ T\in \mathrm{GL}(2,\R) $ on the vector field $ X $ given by the pullback,
			\[ T\cdot X = T X(T^{-1} x), \] 
			induces a right action of $ S \in \mathrm{PGL}(2,\R) $, the group of real M\"{o}bius transformations, on the sets of characteristic directions $ \hat{\C}^{n+1} $ by,
			\[ \hat{\omega}\cdot S = S^{-1}(\hat{\omega}) := (S^{-1}(\omega_1),\dots,S^{-1}(\omega_{n+1})). \]
		\end{proposition}
		\begin{proof}
			We provide the action on $ v = \displaystyle\frac{y}{x} $ and remark that this induces and action on $ \omega $. A coordinate transformation 
			$ T = \left(\begin{array}{cc}
				a & b \\
				c & d \\
			\end{array} \right) $
			induces a transformation $ S $ of $ v $ of the form,
			$\displaystyle S(v) = \frac{c + d v}{a + b v}, $
			clearly in $ \mathrm{PGL}(2,\R) $. A calculation shows that the time derivative of $ v $ is transformed to,
			\begin{equation}
				\dot{v} = \frac{x^{n-1}(d-b v)^{n+1}}{\det T} p\circ S^{-1}(v).
			\end{equation}
			Hence, the new characteristic orbits are given by the roots of the polynomial $ p \circ S^{-1} (v) $. That is, $ T $ transforms the set of roots $ 
			\hat{\omega} = (\omega_1,\dots,\omega_{n+1}) $ to $ S^{-1}(\omega):= (S^{-1}(\omega_1),\dots,S^{-1}(\omega_{n+1})) $.
		\end{proof}
	
		With this proposition, the problem of describing a normalisation is reduced to defining representatives of the orbits $ \mathcal{O}_{\hat{\omega}} $ produced by the induced action. 
		
		Recall that the $ \omega_i  $ are roots of a real degree $ n $ polynomial. Hence, if a root is complex, then the conjugate must also be a root, which further implies that if the order of $ X $, $ n $, is even (odd) then there can only be an even (odd) number of real roots. It is a well-known that there exists a real M\"{o}bius transformation mapping any line or circle with a centre on $ \R $ to any line or circle with centre on $ \R $. Using this, one can deduce the structure of the orbit $ \mathcal{O}_{\hat{\omega}} $ and choose a suitable representative of the orbits to stand as the normalised vector field. For arbitrary $ n $, it is difficult to argue that a choice is undeniably the right one for all possible vector fields. Instead, the choice will depend on the dynamics desired to be emphasised or on any symmetries of the system. For $n=2$ such a choice can be made:
		
		\begin{proposition}\label{prop:normal}
			Suppose $ \hat{\omega} \in \hat{\C}^3 $ are the roots of a degree three real polynomial $ p(x) $ with distinct roots. Then, there exists $ S \in \mathrm{PGL}(2,\R) $ such that $ \hat{\omega} \cdot S $ is of the form,
			\begin{enumerate}
				\item  $ (0,-i,i) $ if $ p $ has one real and two purely imaginary roots.
				\item  $ (0,\infty,1) $ if $ p $ has three real roots.
			\end{enumerate}
			That is, any $ \hat{\omega} \in \hat{\C}^3 $ arising as the roots of a cubic real polynomial is in the orbit of $ (0,-i,i) $ or $ (0,\infty,1) $
		\end{proposition}

%% file: QuadraticVF.tex
\subsection{Regularity of Quadratic Vector Fields}	
		Using the method from Section \ref{sec:Normalisation} a canonical form for the case when $ P,Q $ are homogeneous of degree 2 can be derived. The possible canonical forms are described in the following proposition.		
		\begin{proposition}\label{prop:C0quadratic}
			Let $ X \in \mathcal{V} $ be a homogeneous quadratic vector field. Then $ X $ is $ C^0 $-regularisable if and only if it takes the canonical form
			\begin{equation}\label{eqn:C0quadvf}
				X = \left(\kappa_1 x^2 + \kappa_2 x y - y^2,\ (\kappa_1 + 1) x y + \kappa_2 y^2 \right)
			\end{equation}
			where $ \kappa_1 < 0  $ and $ \kappa_2 \in \R $.
		\end{proposition}
		\begin{proof}
			From Theorem \ref{thm:C0RegularisationHomX} we know that $ X $ is $ C^0 $-regularisable if and only if there is precisely one real root of the directional polynomial $ p(\hat{y}) $. Proposition \ref{prop:normal} asserts that under a general real linear coordinate transformation of $ (x,y) $, the three roots of $ p(\hat{y}) $ can be transformed to $ (0,-i,i) $. That is, $ p(\hat{y}) $ transforms to the canonical form $ \tilde{p}(\hat{y}) = \hat{y}(1+\hat{y}^2) $. Explicitly, let $ X $ take the general form
			\[ X = \left(a_{20} x^2+a_{11} xy + a_{02}y^2,\ b_{20} x^2+b_{11} xy + b_{02}y^2\right).  \]
			Then $ p(\hat y) = -a_{20} \hat{y}^3 + (b_{02}-a_{11})\hat{y}^2 + (b_{11}-a_{20}) \hat{y} +b_{20}. $ So, if $ p(\hat{y}) $ is to take the canonical form, then it must be that $ a_{20} = -1, b_{02}=a_{11},b_{11}=a_{20}+1,b_{20}=0 $. Renaming $ a_{20}=\kappa_1, a_{11}=\kappa_2 $ produces the desired canonical form.
%			\[ \tilde{X} = \left(\kappa_1 x^2 + \kappa_2 x y - y^2,\ (\kappa_1 + 1) x y + \kappa_2 y^2\right). \]
%			
			To finish the proof, note from Theorem \ref{thm:C0RegularisationHomX} $ X $ has $ C^0 $ regularity if and only if the characteristic value $ r^* < 0 $. By computing the blow-up, it is easily seen that $ \kappa_1 = r^* $, hence $ \kappa_1<0 $. In particular, when $ \kappa_1 > 0 $ the block map is not even continuous.
		\end{proof}
		
\begin{figure}
			\begin{tabular}{cc}
				\includegraphics[width=0.48\textwidth]{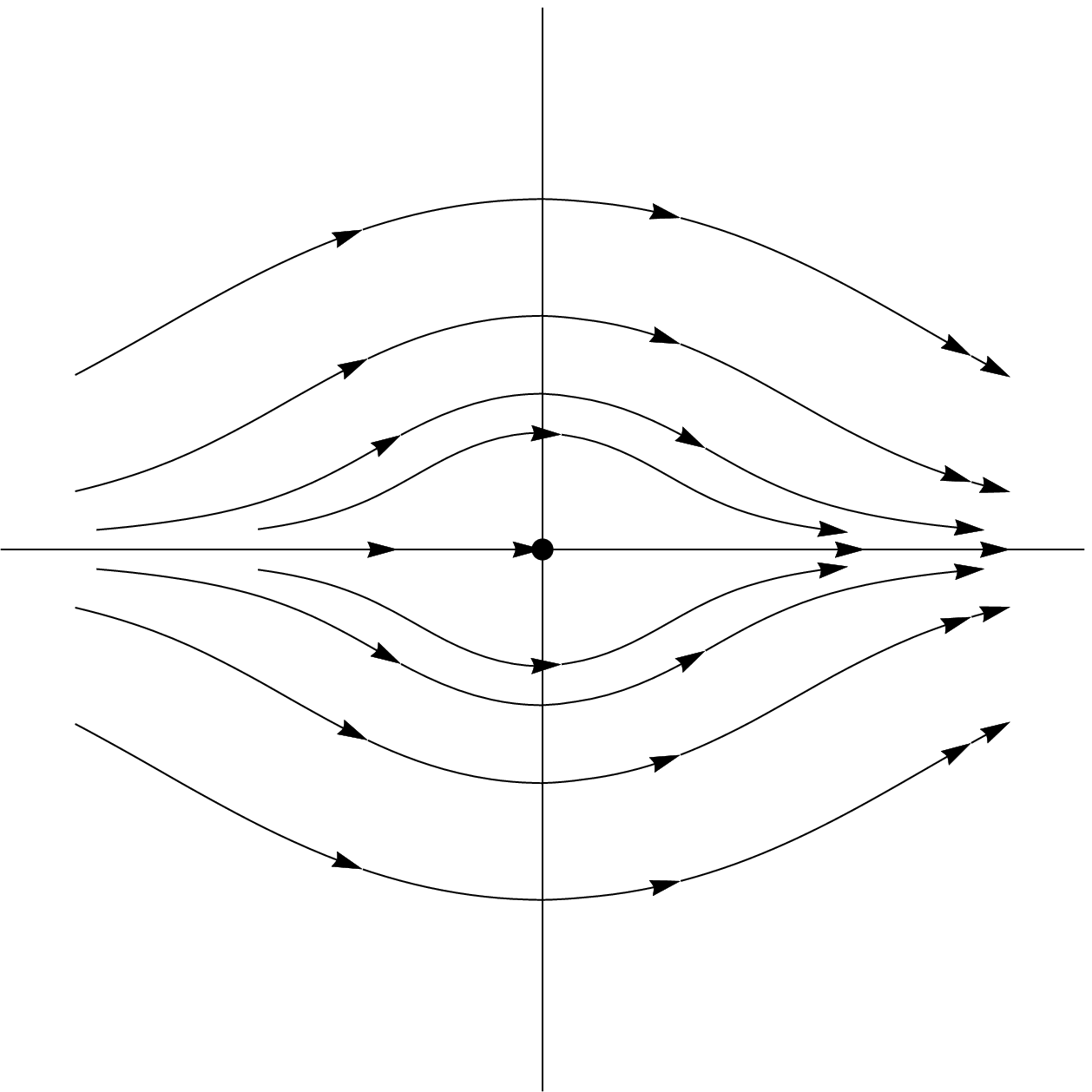} &   \includegraphics[width=0.48\textwidth]{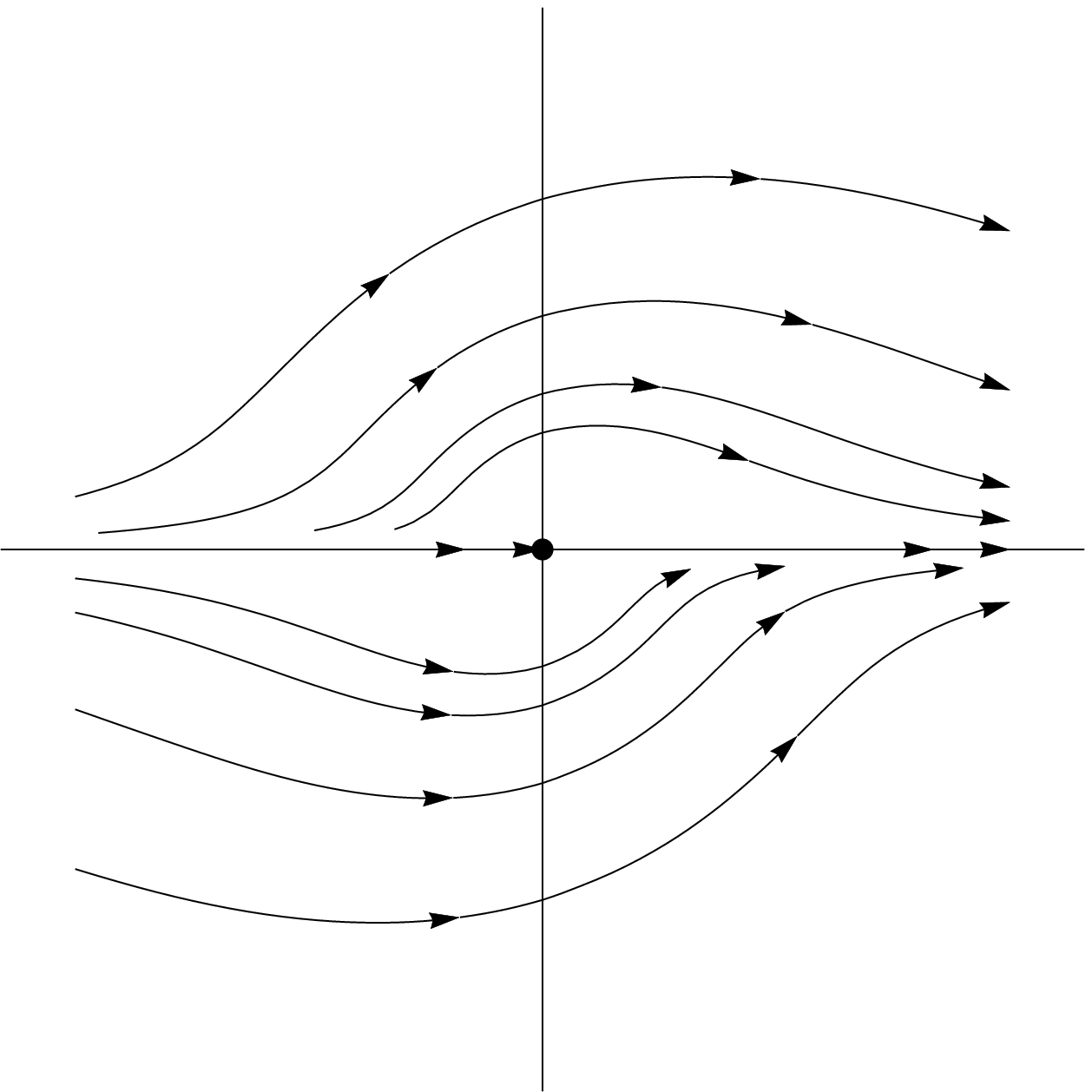} \\
				(a) $ C^\infty : \left( \frac{1}{3}x^2 + y^2, -\frac{2}{3} x y \right) $ & (b) $ C^{0,1} : \left( \frac{1}{3}x^2 + y^2 + \frac{1}{2} x y, -\frac{2}{3} x y+\frac{1}{2} y^2 \right) $ \\[6pt]
				\includegraphics[width=0.48\textwidth]{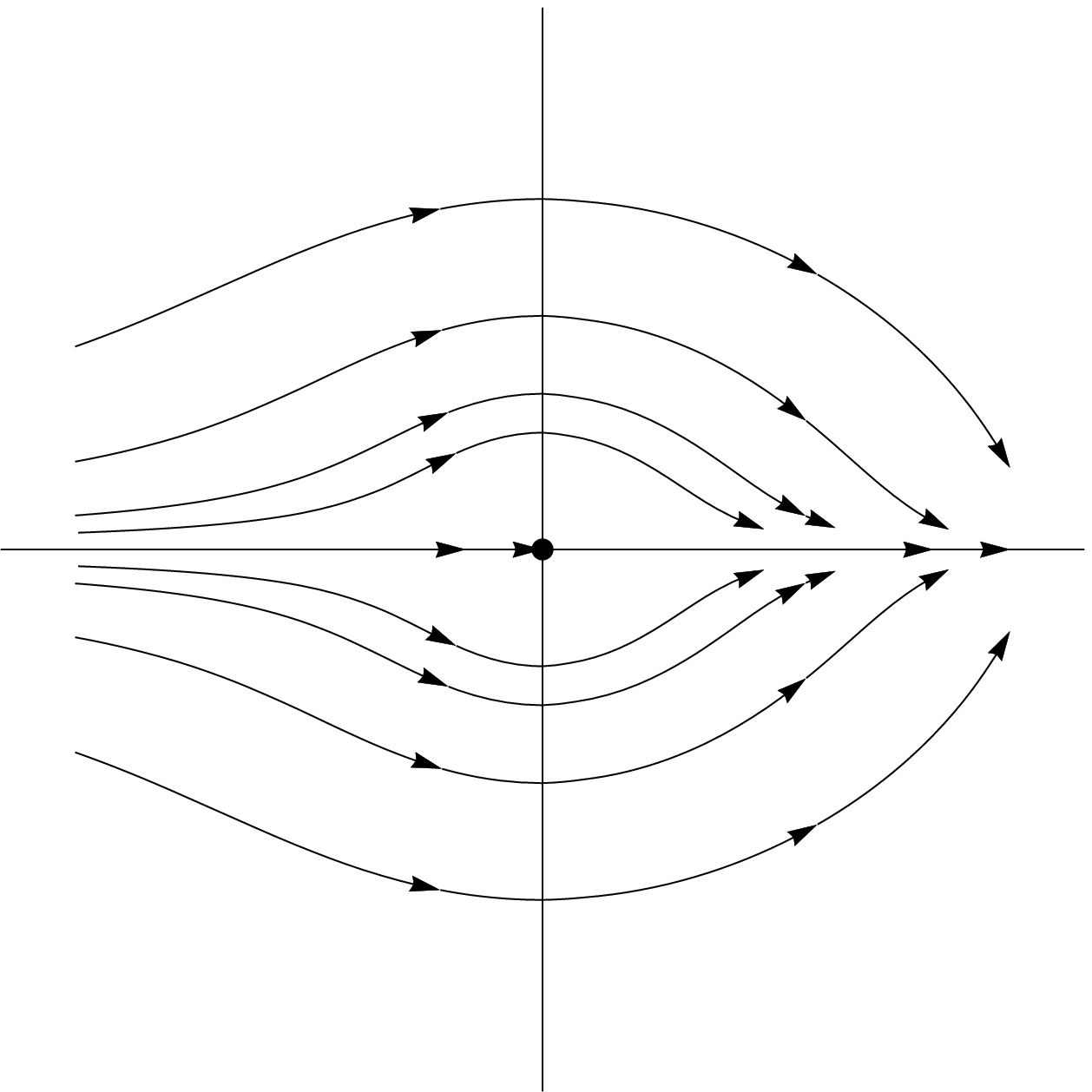} &   \includegraphics[width=0.48\textwidth]{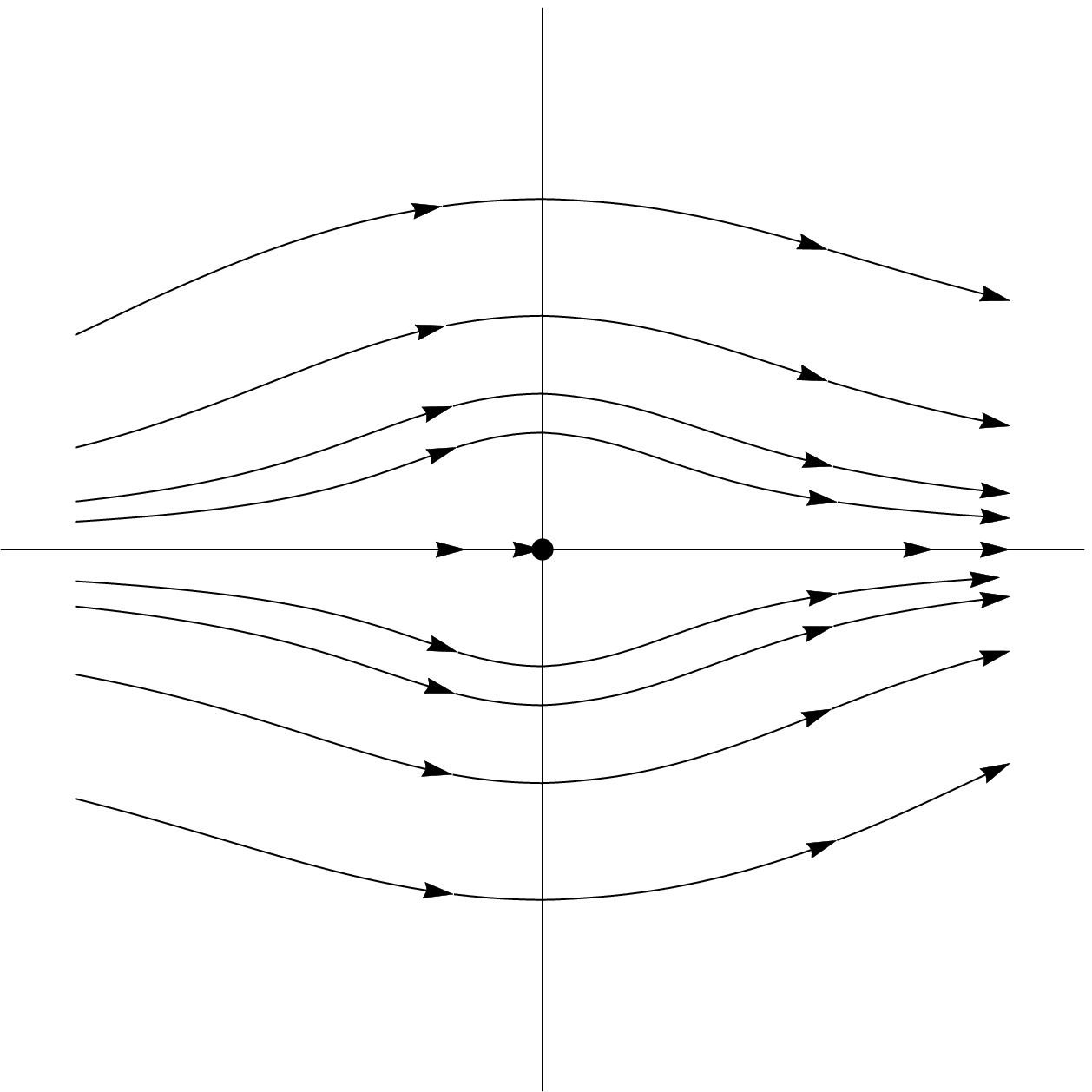} \\
				(c) $ C^{4/3} : \left( \frac{1}{3}x^2 + y^2+\frac{1}{2} x^3, -\frac{2}{3} x y \right) $ & (d) $ C^{1,a} : \left( \frac{1}{2}x^2 + y^2+\frac{1}{2} x^3 y, -\frac{1}{2} x y \right) $ \\[6pt]
			\end{tabular}
			\caption{Examples of the mechanisms that create finite regularity of the block map. Each corresponds to an example in Section \ref{sec:examples}. Figure (a) has a smooth block map. Figure (b) corresponds to the upper and lower trajectories disagreeing at the first derivative. Figure (c) is an example of when the order of transition $ k < (m-1)q $. Figure (d) shows a vector field which is dominated by the order of resonance $ k > (m-1)q $.}\label{fig:examples}
		\end{figure}

		Using the results from Section \ref{sec:RegularisingMap} we can determine which of the $ C^0 $-regularisable $ X \in \mathcal{V} $ pertain higher order regularity.
		\begin{proposition}\label{prop:quadraticSmooth}
			Let $ X \in \mathcal{V} $ be a homogeneous quadratic vector field. Then $ X $ is $ C^\infty $-regularisable if and only if it takes the canonical form,
			\begin{equation}\label{eqn:C1form}
				X = \left( \kappa x^2 - y^2,\ (\kappa + 1) x y \right)
			\end{equation}
			where $ \kappa < 0  $ is the unique characteristic value of $ X $. If $ \kappa_2 \not = 0 $ then $ X $ is $ C^{0,1} $-regularisable.
		\end{proposition}	
		\begin{proof}
			From the previous Proposition \ref{prop:C0quadratic} we already know the form of the $ C^0 $-regularisable quadratic homogeneous vector fields, namely \eqref{eqn:C0quadvf}. Now, for $ C^1 $ regularity, Corollary \ref{cor:C0regularisation} asserts that the Cauchy Principal value of the integral \eqref{eqn:C1reg} must be 0. It can be seen that $ q_2(x) = (-1-x^2) $ and  $ Q_s(x,1) = (1+\kappa_1)x + \kappa_2 $ . Then, either by directly calculating the integral and taking limits or using the residue theorem, it is found that,
			\begin{equation}
				\int_{-\infty}^{\infty} -\frac{(1+\kappa_1)x + \kappa_2}{1+x^2} dx = -\kappa_2 \pi.
			\end{equation}
			Hence, for $ C^1 $ regularity of $ X $ it must be that $ \kappa_2 = 0 $. 
			When $\kappa_2 \not = 0$ the block map is only $C^{0,1}$.
			Finally, if $ \kappa_2 = 0 $ then $ X $ has a discrete reversing symmetry $ (x,t) \mapsto (-x,-t) $. By choosing $ \Sigma_0 = \set{x=-c}, \Sigma_3 = \set{x=c} $ a regularising map $ \bar{\pi} $ is produced that is the identity and the result is concluded.
		\end{proof}
		Figure \ref{fig:examples} shows several examples of regularisable vector fields. In particular, when $ \kappa_1 = 1/3 $ and $ \kappa_2 = 0 $ we have in Figure \ref{fig:examples} (a) an example of a $ C^\infty $-regularisable singularity. In (b) we have $ \kappa_2 \neq 0 $. Notice how above the $ x $-axis the orbits are expanding. Where as below the $ x $-axis there is contraction. In this case the block-map is only Lipschitz.
	
	\subsection{Examples of Finite Regularity}\label{sec:examples}
	
		\begin{example}[Perturbed Toy]\label{exam:PerturbedToy}
			It has been shown in earlier sections that the toy example \eqref{exam:Toy} is $ C^\infty $-regularisable with characteristic value $ r^* = -1/3 $. From Proposition \ref{prop:quadraticSmooth} the canonical form of the vector field is,
			\[ X = \left( -\frac{1}{3} x^2 - y^2,\ \frac{2}{3} x y \right). \]
			We show that a perturbation of the toy,
			
			\[ X_{p} = \left(-\frac{1}{3} x^2 - y^2 + \lambda x^3,\ \frac{2}{3} x y \right) \]
			is precisely $ C^{1+1/3} $-regularisable for $ \lambda \neq 0 $. A plot of the vector field is given in Figure \ref{fig:examples}.
			
			The steps to the proof are outlined in Section \ref{ssec:iterativemethod}. In order to determine the $ C^{1+1/3} $ regularity we must check that there are no resonances of order $ m \leq 2 $ and $ f $ has an order of transition $ k = 2 $. Hence, the normal transformation $ \phi:\hat{X}\to X_N $ must be computed to degree 5. The Dulac maps $ d_1,d_2 $  will be computed from the normal form and then, by solving the $ 2^{nd} $ order variational equations, we will obtain the smooth transition map $ F $ to order $ 2 $ and show the order of the transition is $ 2 $ as required.
			
			The $ x $-directional blowup $ \hat{X}_p $, is computed from the transformation $ x=\hat{x}, y=\hat{x}\hat{y} $ followed by the time rescaling $ d\tau=\hat{x} dt $. This results in
			\begin{equation}
				\hat{X}_p = \left(-\frac{1}{3}\hat{x}(1 + 3\hat{y}^2)+\lambda \hat{x}^2,\ \hat{y}+\hat{y}^3 -\lambda \hat{x} \hat{y}\right).			
			\end{equation}
			
			The normal form transformation $ \phi^{-1}:(w,v)\mapsto (\hat{x},\hat{y}) $ up to degree 5 is computed as
			\begin{equation}
				\begin{aligned}
					\hat{x} &= w\left(1-\frac{1}{2} v^2-\frac{1}{8} v^4\right) - \lambda w^2\left(3 -\frac{12}{5} v^2\right) + \lambda^2 w^3 \left( 9 -\frac{27}{20} v^2 \right)-27\lambda^3 w^4 +81 \lambda^4 w^5 + \dots \\
					\hat{y} &= v+\frac{1}{2}v^2+\frac{3}{8}v^5+\lambda w \left( 3v +\frac{27}{5} v^3 \right) + \frac{81}{4}\lambda^2 v^3 w^2 + \dots	
				\end{aligned}
			\end{equation}
			Under this transformation $ \hat{X} $ is transformed into the normal form
			\begin{equation}
				X_N = (-1/3 w,v),
			\end{equation}
			and Dulac maps $ d_{1,\epsilon}^+:\tau_0 \to \sigma_1^\epsilon $, $ d_{2,\epsilon}^+: \sigma_2^\epsilon \to \tau_3 $ are simply,
			\begin{equation}
				d_{1,\epsilon}^+(w_0) = \epsilon^{-1/3} w_0^{1/3} + \dots,\quad d_{2,\epsilon}^+(v_2) = \epsilon^{-1} v_2^{3} +\dots .
			\end{equation}
			As no resonance terms have appeared to this degree 5, we know the order of resonance $ m > 2 $.
			
			Now, we wish to compute the order of transition $ k $. To compute the Taylor series of the smooth transition map $ F_\epsilon:\overline{\Sigma}^\epsilon_1\to\overline{\Sigma}^\epsilon_2 $ we can compute the derivatives using higher order variational equations. From Example \ref{exam:toyC1}, the first derivative of $ F_\epsilon $ is $ F'_\epsilon(0) = 1 $. The higher order variational equations become difficult to solve as the order increases, yet they are always linear. In order to give the form for the 2nd order variational equation, we first compute $ \bar{X} $:
			\begin{equation}
				\bar{X} = \left( -(1+\bar{x}^2) + \lambda \bar{y}\bar{x}^3,\ \frac{2}{3} \bar{x}\bar{y} \right).
			\end{equation}
			Then $ \displaystyle \frac{d\bar{y}}{d\bar{x}} $ is given by
			\[ \frac{d\bar{y}}{d\bar{x}} = -\frac{2}{3} \frac{\bar{x}}{(1+\bar{x}^2)}\bar{y} - \frac{2}{3} \lambda \frac{\bar{x}^4}{(1+\bar{x}^2)^2} \bar{y}^2 + O(\bar{y}^3).  \]
			For the variational equations, assume a solution of the form 
			\begin{equation}\label{eqn:vareqn}
				\bar{y}(\bar{x},\bar{y}_0) = g_1(\bar{x}) \bar{y}_0 + \frac{1}{2} g_2(\bar{x}) \bar{y}_0^2 + O(\bar{y}_0)^2,
			\end{equation} 
			$ g_1(0)=1 $ and $ g_i(0) = 0 $ for $ i >1 $. The variational equations give the coefficients as,
			\begin{equation} \label{eqn:diffeqns}
				\begin{aligned}
					g_1(\bar{x}) &= \exp\left(-\int_0^{\bar{x}} \frac{2}{3} \frac{x}{(1+x^2)} dx\right) \\
					g_2(\bar{x}) &= g_1(\bar{x})\int_{0}^{\bar{x}} \frac{2}{3} \lambda \frac{x^4}{(1+x^2)^2} g_1(x) dx
				\end{aligned}
			\end{equation}
			Now, we want to compute the map between $ \bar{\Sigma}_1^\epsilon,\bar{\Sigma}_2^\epsilon $ the images of $ \sigma_i^\epsilon $ in $ \bar{X} $. Using $ \phi, P_x,P_y $ we obtain these sections as
			\begin{equation}\label{eqn:secs}
				\begin{aligned}
					\bar{\Sigma}_1^\epsilon &= \set{\left( \frac{1}{\epsilon} - \frac{3\lambda}{\epsilon} w_1 + O(w_1^2),\ \epsilon w_1 + O(w_1^3) \right):\  w_1\in[0,1)},\\
					\bar{\Sigma}_2^\epsilon &= \set{\left( -\frac{1}{\epsilon} - \frac{3\lambda}{\epsilon} w_2 + O(w_2^2),\ \epsilon w_2 + O(w_2^3) \right):\ w_2\in[0,1)} \,.
				\end{aligned}
			\end{equation}
			Using equation \eqref{eqn:vareqn} twice, once with $ (\bar{x}(w_1),\bar{y}(w_1)) $ on $ \bar{\Sigma}_1^\epsilon $ and a second time with $ (\bar{x}(w_2),\bar{y}(w_2)) $ on $ \bar{\Sigma}_2^\epsilon $, we have an implicit relation between $ w_1 $ and $ w_2 $. Exploiting this and solving the integrals \eqref{eqn:diffeqns} we obtain the form of $ w_2(w_1) = F_\epsilon(w_1) $,
			\begin{equation}\label{eqn:Feps}
				F_\epsilon(w_1) = (1+O(\epsilon^2))w_1 -\left(\frac{9\lambda \sqrt{\pi}\Gamma(-1/6)}{8\Gamma(1/3)} \epsilon^{1/3} + O(\epsilon^2)\right) w_1^2 + O(w_1^3)
			\end{equation}
			The details of this calculation are found in Appendix \ref{sec:calculation}.
			
			Finally, composing $ d_{2,\epsilon}\circ F_\epsilon\circ d_{1,\epsilon}(x) $ and taking $ \epsilon \to 0 $ we obtain the transition,
			\begin{equation}\label{eqn:blockmap}
				\pib^+(x) = x - \frac{27 \sqrt{\pi} \Gamma{(-1/6)}}{8 \Gamma{(1/3)}} x^{1+1/3}+\dots
			\end{equation} 
			and conclude that the perturbed toy is precisely $ C^{4/3} $-regularisable.
		\end{example}
	
		\begin{example}[Resonant Toy]
			We show that the vector field
			\[ X_{R} = \left(\frac{1}{2} x^2 + y^2 + \lambda x^3y,\ \frac{1}{2} x y \right) \]
			has a block map whose first derivative is Log-Lipschitz, that is, provided $ \lambda \neq 0 $, $ X_R $ is $ C^{1,a} $-regularisable for $ 0 \leq a < 1 $. 
			
			This is an example of when there is a finite order of resonance $ m $ that dominates the block map. A plot of the vector field is given in Figure \ref{fig:examples}. 
			
			Firstly, the leading order homogeneous component of $ X_R $ is of the form \eqref{eqn:C1form} and so it is concluded that it is at least $ C^0 $ regularisable. Note that the characteristic value $ r^* = 1/2 $, hence $ p = 1, q = 2 $.
			
			To determine the $ C^{1,a} $ regularity we must check that there is a resonance of order $ m = 2 $ and $ f $ has an order of transition $ k \geq (m-1)q = 2 $. Therefore the normal transformation $ \phi:\hat{X}\to X_N $ must be computed to degree 5.
			
			The $ x $-directional blowup $ \hat{X}_R $ is
			\begin{equation}
			\hat{X}_R = \left(\frac{1}{2}\hat{x}(1 + \hat{y}^2)+\lambda \hat{x}^3 \hat{y},\ -\hat{y}(1+\hat{y}^2) -\lambda \hat{x}^2 \hat{y}^2\right).			
			\end{equation}
			
			The normal form transformation $ \phi^{-1}:(w,v)\mapsto (\hat{x},\hat{y}) $ up to degree 5 is computed as
			\begin{equation}
				\begin{aligned}
					\hat{x} &=  w -\frac{1}{2}w v^2 - \frac{1}{8}w v^4 \\
					\hat{y} &=  v+\frac{1}{2} v^3 +\frac{3}{8} v^5
				\end{aligned}
			\end{equation}
			Under this transformation $ \hat{X} $ is transformed into the normal form
			\begin{equation}
				X_N = \left( w+ \lambda w^3 v, -v - \lambda w^2 v^2 \right).
			\end{equation}
			As described in the proof of Proposition \ref{prop:normalform}, the order of resonance is computed by introducing the variable $ u = w^1 v^2 $ yielding
			\begin{equation}
				\dot{u} = -a u^2 + O(u^3).
			\end{equation}
			We can conclude that the order of resonance is $ m = 2 $ as desired. Further, the coefficient $ \alpha_{1,2} = -a $ and it can be confirmed in a similar fashion that $ \alpha_{2,2} = -2a $.
			
			Lastly, we must check that $ k \geq 3 $. The $ y $-directional blow up is given by
			\begin{equation}
			\bar{X} = \left( (1+\bar{x}^2) + \lambda \bar{y}^2\bar{x}^3,\ -\frac{1}{2} \bar{x}\bar{y} \right)
			\end{equation}
			and so $ \displaystyle \frac{d\bar{y}}{d\bar{x}} $ is
			\begin{equation}
				\frac{d\bar{y}}{d\bar{x}} = -\frac{1}{2} \frac{\bar{x}}{(1+\bar{x}^2)}\bar{y} + \frac{1}{2} \lambda \frac{\bar{x}^4}{(1+\bar{x}^2)^2} \bar{y}^3 + O(\bar{y}^4). 
			\end{equation}
			As the coefficient of $ \hat{y}^2 $ is 0, the solution to the second order variational equation will be $ 0 $. Hence, we can conclude that the order of transition $ k \geq 3 $.
			
			Using Lemma \ref{lem:leadingorder} the upper block map is
			\begin{equation}
				\pib^+(x) = x + \lambda x^2 \ln x + \dots
			\end{equation} 
			This confirms that $ X_R $ is $ C^{1,a} $-regularisable as desired.
		\end{example}

%% file: Appendix.tex
\section{Calculation of Order of Transition for Perturbed Toy}\label{sec:calculation}

This appendix contains more details on the final calculations of Example \ref{exam:PerturbedToy}. The coefficients of the transition map $ F_\epsilon(w_1) $ in \eqref{eqn:Feps} and the final upper block map $ \pib^+(x) $ in \eqref{eqn:blockmap} are calculated.

We begin by first evaluating the integrals \eqref{eqn:diffeqns} to get the solutions to the first and second order variational equations
\begin{equation} 
	\begin{aligned}
		g_1(\bar{x}) &= \exp\left(-\int_0^{\bar{x}} \frac{2}{3} \frac{x}{(1+x^2)} dx\right)\\
			&= \left(1+\bar{x}^2\right)^{-1/3}\\
		g_2(\bar{x}) &= g_1(\bar{x})\int_{0}^{\bar{x}} \frac{2}{3} \lambda \frac{x^4}{(1+x^2)^2} g_1(x) dx\\
		&= -\frac{1}{4}\lambda \bar{x}\left( 9\left(1+\bar{x}^2\right)^{-1/3} {}_2F_1\left(1,\frac{7}{6};\frac{3}{2};-\bar{x}^2\right)-\left(1+\bar{x}^2\right)^{-5/3}\left(9+11\bar{x}^2\right) \right).
	\end{aligned}
\end{equation}
We seek to compute $ F_\epsilon $, the transition between $ \bar{\Sigma}_1 $ and $ \bar{\Sigma}_2 $. These sections are parameterised as given in \eqref{eqn:secs} and we write this parametrisation as $ (\bar{x}_1(w_1),\bar{y}_1(w_1)) $ and $ (\bar{x}_2(w_2),\bar{y}_2(w_2)) $. Now, using \eqref{eqn:vareqn}, we have the two relations,
\begin{equation}\label{eqn:us}
	\begin{aligned}
		\bar{y}_1(w_1) &= g_1(\bar{x}_1(w_1)) \bar{y}_0 + \frac{1}{2}g_2(\bar{x}_1(w_1)) \bar{y}_0^2 + O(\bar{y}_0)^2\\
		\bar{y}_2(w_2) &= g_1(\bar{x}_2(w_2)) \bar{y}_0 + \frac{1}{2}g_2(\bar{x}_2(w_2)) \bar{y}_0^2 + O(\bar{y}_0)^2.
	\end{aligned}
\end{equation}
The map $ F(w_1) $ is exactly $ w_2(w_1) $ which is implicitly defined through these equations \eqref{eqn:us}. In order to get the Taylor series about $ w_1 = 0 $, the derivatives are computed by considering $ w_2 = w_2(w_1) $ and $ y_0 = y_0(w_1) $, differentiating both relations in \eqref{eqn:us} with respect to $ w_1 $, evaluating at $ w_1 = y_0 =0 $ and solving the resulting system of equations in $ y_0'(0), w_2'(0) $ for $ w_2'(0) $. Similarly, taking a second derivative we can compute $ w_2''(0) $
\begin{align}
	w_2'(0)		&= 1+O(\epsilon^2) \\
	w_2''(0)	&= \frac{9}{2}\lambda\left( 3-\frac{1}{\epsilon^2}{}_2F_1\left( 1,\frac{7}{6}; \frac{3}{2}; -\frac{1}{\epsilon^2} \right) \right)+ O(\epsilon^2).
\end{align}
To complete the derivation, the asymptotic series at $ \epsilon=0 $ of the hypergeometric function $ {}_2F_1\left( 1,\frac{7}{6}; \frac{3}{2}; -\frac{1}{\epsilon^2} \right) $ must be computed. To do this, we make use of the series definition and some identities,
\begin{align}
	{}_2F_1(a,b;c;z) &= \sum_{s=0}^{\infty}\frac{(a)_s (b)_s}{(c)_s s!} z^s,\quad |z|<1 \\
	{}_2F_1(a,b;c;-z) &= \frac{\Gamma(b-a)\Gamma(c)}{\Gamma(b)\Gamma(c-a)} z^{-a} {}_2F_1\left(a,a-c+1;a-b+1;\frac{1}{z}\right) \\
	&\qquad+ \frac{\Gamma(a-b)\Gamma(c)}{\Gamma(a)\Gamma(c-b)} z^{-b} {}_2F_1\left(b,b-c+1;b-a+1;\frac{1}{z}\right),\nonumber
\end{align}
where $ (x)_s $ is the Pochhammer symbol. Using the second relation, then the first, it can be seen that
\begin{equation}
	{}_2F_1\left( 1,\frac{7}{6}; \frac{3}{2}; -\frac{1}{\epsilon^2} \right) = 3 \epsilon^2 + \frac{\sqrt{\pi}\Gamma(-1/6)}{2 \Gamma(1/3)}\epsilon^{7/3} + O(\epsilon^4).
\end{equation}
Substituting into the expression for $ w_2''(0) $ we have
\begin{equation}
	\begin{aligned}
		F_\epsilon(w_1)	&= w_2'(0) w_1 + \frac{1}{2} w_2''(0) w_1^2 + O(w_1^3)\\
						&= (1+O(\epsilon^2))w_1 -\left(\frac{9\lambda \sqrt{\pi}\Gamma(-1/6)}{8\Gamma(1/3)} \epsilon^{1/3} + O(\epsilon^2)\right) w_1^2 +O(w_1^3),
	\end{aligned}
\end{equation}
as desired.

Lastly, to compute the upper block map $ \pib^+ $, the Dulac maps $ d_{1,\epsilon},d_{2,\epsilon} $ and the smooth transition $ F_\epsilon $ must be composed and the limit $ \epsilon \to 0 $ taken. That is,
\begin{align*}
	\pib^+	&= \lim_{\epsilon\to 0} d_{2,\epsilon}\circ F_\epsilon \circ d_{1,\epsilon}(x) \\
			&= \lim_{\epsilon\to 0}  d_{2,\epsilon}\circ F_\epsilon(\epsilon^{-1/3} x^{1/3} + \dots) \\
			&= \lim_{\epsilon\to 0}  d_{2,\epsilon}\left( (1+O(\epsilon^2))\epsilon^{-1/3} x^{1/3} -\left(\frac{9\lambda \sqrt{\pi}\Gamma(-1/6)}{8\Gamma(1/3)} \epsilon^{1/3} + O(\epsilon^2)\right) \epsilon^{-2/3} x^{2/3} + \dots \right) \\
			&= \lim_{\epsilon\to 0}  \epsilon^{-1}\left( (1+O(\epsilon^2))\epsilon^{-1/3} x^{1/3} -\left(\frac{9\lambda \sqrt{\pi}\Gamma(-1/6)}{8\Gamma(1/3)} \epsilon^{1/3} + O(\epsilon^2)\right) \epsilon^{-2/3} x^{2/3} + \dots \right)^3 \\
			&= \lim_{\epsilon\to 0} x\left( (1+O(\epsilon^2)) -\left(\frac{27\lambda \sqrt{\pi}\Gamma(-1/6)}{8\Gamma(1/3)} + O(\epsilon^{5/3})\right)  x^{1/3} + \dots \right) \\
			&= x\left( 1 -\frac{27\lambda \sqrt{\pi}\Gamma(-1/6)}{8\Gamma(1/3)} x^{1/3} + \dots  \right).
\end{align*}

%% file: Deg2VF.bbl
\begin{thebibliography}{10}

\bibitem{andronov1974qualitative}
AA~Andronov, EA~Leontovich, II~Gordon, AG~Maier, and Martin~C Gutzwiller.
\newblock Qualitative theory of second-order dynamic systems.
\newblock {\em Physics Today}, 27:53, 1974.

\bibitem{bakker2011existence}
Lennard~F Bakker, Tiancheng Ouyang, Duokui Yan, and Skyler Simmons.
\newblock Existence and stability of symmetric periodic simultaneous binary
  collision orbits in the planar pairwise symmetric four-body problem.
\newblock {\em Celestial Mechanics and Dynamical Astronomy}, 110(3):271--290,
  2011.

\bibitem{Belitskii2002}
G.~Belitskii.
\newblock {$C^\infty$}-normal forms of local vector fields.
\newblock {\em Acta Appl. Math.}, 70(1-3):23--41, 2002.
\newblock Symmetry and perturbation theory.

\bibitem{Brunella1990}
Marco Brunella and Massimo Miari.
\newblock Topological equivalence of a plane vector field with its principal
  part defined through {N}ewton polyhedra.
\newblock {\em J. Differential Equations}, 85(2):338--366, 1990.

\bibitem{Bruno1989}
Alexander~D. Bruno.
\newblock {\em Local methods in nonlinear differential equations}.
\newblock Springer-Verlag, Berlin, 1989.

\bibitem{CG17}
G.~Calsamiglia and Y.~Genzmer.
\newblock Classification of regular dicritical foliations.
\newblock {\em Ergod. Th. Dynam. Sys.}, 37:1443--1479, 2017.

\bibitem{conley1971isolated}
Charles Conley and Robert Easton.
\newblock Isolated invariant sets and isolating blocks.
\newblock {\em Transactions of the American Mathematical Society},
  158(1):35--61, 1971.

\bibitem{Dumortier1977}
Freddy Dumortier.
\newblock Singularities of vector fields on the plane.
\newblock {\em J. Differential Equations}, 23(1):53--106, 1977.

\bibitem{easton1971regularization}
Robert Easton.
\newblock Regularization of vector fields by surgery.
\newblock {\em Journal of Differential Equations}, 10(1):92--99, 1971.

\bibitem{easton1972topology}
Robert Easton.
\newblock The topology of the regularized integral surfaces of the 3-body
  problem.
\newblock {\em Journal of Differential Equations}, 12(2):361--384, 1972.

\bibitem{elbialy1990collision}
Mohamed~Sami ElBialy.
\newblock Collision singularities in celestial mechanics.
\newblock {\em SIAM journal on mathematical analysis}, 21(6):1563--1593, 1990.

\bibitem{elbialy1993simultaneous}
Mohamed~Sami Elbialy.
\newblock Simultaneous binary collisions in the collinear n-body problem.
\newblock {\em Journal of differential equations}, 102(2):209--235, 1993.

\bibitem{Ilyashenko2008}
Yulij Ilyashenko and Sergei Yakovenko.
\newblock {\em Lectures on analytic differential equations}, volume~86 of {\em
  Graduate Studies in Mathematics}.
\newblock American Mathematical Society, Providence, RI, 2008.

\bibitem{Jaurez-Rosas2017}
Jessica~Ang\'elica Jaurez-Rosas.
\newblock Real-formal orbital rigidity for germs of real analytic vector fields
  on the real plane.
\newblock {\em J. Dyn. Control Syst.}, 23(1):89--109, mar 2017.

\bibitem{manosa2002center}
Victor Ma{\~n}osa.
\newblock On the center problem for degenerate singular points of planar vector
  fields.
\newblock {\em International Journal of Bifurcation and Chaos},
  12(04):687--707, 2002.

\bibitem{martinez1999simultaneous}
Regina Mart{\'\i}nez and Carles Sim{\'o}.
\newblock Simultaneous binary collisions in the planar four-body problem.
\newblock {\em Nonlinearity}, 12(4):903, 1999.

\bibitem{martinez2000degree}
Regina Mart{\'\i}nez and Carles Sim{\'o}.
\newblock The degree of differentiability of the regularization of simultaneous
  binary collisions in some n-body problems.
\newblock {\em Nonlinearity}, 13(6):2107, 2000.

\bibitem{mcgehee1974triple}
Richard McGehee.
\newblock Triple collision in the collinear three-body problem.
\newblock {\em Inventiones mathematicae}, 27(3):191--227, 1974.

\bibitem{OBRGV14}
L.~Ortiz-Bobadilla, E.~Rosales-Gonzales, and S.M. Voronin.
\newblock Formal and analytic normal forms of germs of holomorphic nondicritic
  foliations.
\newblock {\em J. of Singularities}, 9:168--192, 2014.

\bibitem{punovsevac2012regularization}
Predrag Puno{\v{s}}evac and Qiudong Wang.
\newblock Regularization of simultaneous binary collisions in some
  gravitational systems.
\newblock {\em The Rocky Mountain Journal of Mathematics}, pages 257--283,
  2012.

\bibitem{roussarie1995bifurcations}
Robert~H Roussarie.
\newblock {\em Bifurcations of Planar Vector Fields and Hilbert's 16th Problem:
  20o Coloquio Brasileiro de Matematica, IMPA 24-28 Julho, 1995}.
\newblock Instituto de Matematica Pura e Aplicada, 1995.

\bibitem{siegel2012lectures}
Carl~L Siegel and J{\"u}rgen~K Moser.
\newblock {\em Lectures on Celestial Mechanics: Reprint of the 1971 Edition}.
\newblock Springer Science \& Business Media, 2012.

\bibitem{VOBRG10}
S.M. Voronin, L.~Ortiz-Bobadilla, and E.~Rosales-Gonzales.
\newblock Thom's problem for orbital analytic classification of degenerate
  singular points of holomorphic vector fields in the plane.
\newblock {\em Doklady Mathematics}, 82(2):759--761, 2010.

\end{thebibliography}
